\documentclass[a4paper,11pt]{article}
\usepackage[english]{babel}
\usepackage[utf8]{inputenc}
\usepackage[T1]{fontenc}
\usepackage{amsmath,amsfonts,amssymb,amsthm,mathtools}
\usepackage{graphicx}
\usepackage{array}
\usepackage{booktabs}
\usepackage{multirow}
\usepackage{subcaption}
\usepackage{tabularx}
\usepackage{colortbl}
\usepackage{orcidlink}
\usepackage{float}
\usepackage{comment}
\usepackage{tikz}
\usetikzlibrary{arrows.meta,positioning,fit,shapes.geometric}
\usepackage{pgfplots}
\pgfplotsset{compat=1.18}
\usepackage{algorithm}
\usepackage{algpseudocode}

\makeatletter
\renewcommand{\ALG@beginalgorithmic}{\small}
\makeatother

\newtheorem{theorem}{Theorem}[section]

\newtheorem{definition}[theorem]{Definition}

\begin{document}

\title{ADPSO-ERLS: A Hybrid Discrete PSO with Enhanced\\
	Local Search for the Traveling Salesman Problem}
\author{{Anier Soria-Lorente \orcidlink{0000-0003-3488-3094}}$^{1}$,
{Jean-Marie Vilaire \orcidlink{0000-0002-4012-6471}}$^{2}$,\\
{Junior Michel \orcidlink{0009-0005-8450-2757}}$^{1}$, {L\'{a}zaro Alvarado-Barrios \orcidlink{0000-0002-6030-9582}}$^{3}$ \\
\\
$^{1}$Department of Quantitative Methods, Loyola University,\\
Avda. de las Universidades, 2, Dos Hermanas, Seville,\\
41704, Andalusia, Spain\\
asoria@uloyola.es, jmichel@uloyola.es\\
$^{2}$Mobile Computing and Networking Research Laboratory (LARIM),\\
Department of Computer and Software Engineering,\\
Polytechnique Montr\'eal, Montr\'eal, QC H3T 1J4, Canada\\
jean-marie.vilaire@polymtl.ca\\
$^{3}$Department of Engineering, Loyola University,\\
Avda. de las Universidades, 2, Dos Hermanas, Seville,\\
41704, Andalusia, Spain\\
lalvarado@uloyola.es\\
}
\maketitle{}
\begin{abstract}
	The Traveling Salesman Problem is a canonical setting for studying how a
	population-based method should allocate a fixed search budget between
	exploration and progressively stronger local intensification. We propose
	ADPSO-ERLS, a discrete swarm algorithm that treats this allocation as an
	explicit, tunable design variable. It couples memory-guided swap mutation,
	heterogeneous initialization, selective candidate-restricted 2-opt during
	evolution, and an incumbent-only final refinement combining candidate-restricted
	and optional full 2-opt with double-bridge perturbations. The method is
	PSO-inspired, using personal and global memories yet dispensing with velocity,
	inertia, and acceleration coefficients. All six algorithms are implemented in
	Rust, run on identical hardware, and stopped at a strict, recorded limit of
	100,000 candidate-solution assessments, so that programming language, hardware,
	and evaluation budget are held common across methods; wall-clock time is
	reported separately because equal assessment counts need not correspond to
	equal arithmetic work. Over 50 runs on five symmetric TSPLIB instances under the integer
	\texttt{EUC\_2D} convention, ADPSO-ERLS attains the lowest best and mean cost on
	every instance, with best-tour Gap of $1.93$--$4.17\%$ and relative error of
	$3.04$--$5.50\%$. It ranks first under the Friedman test, and all twenty-five
	multiplicity-controlled Wilcoxon comparisons favor it with large, near-complete
	distributional separation. A paired ablation with common seeds links
	initialization, in-run local search, and final refinement to quality gains,
	while candidate restriction chiefly cuts runtime, by up to a factor of roughly
	$38$. Further experiments up to $16{,}862$ cities keep best-tour Gaps below
	$6.7\%$, solving the largest case in under eleven minutes.
	
	\vspace{0.3cm}
	
	
	{\footnotesize Corresponding author: Anier Soria-Lorente}

\end{abstract}



\section{Introduction}\label{sec:intro}
The Traveling Salesman Problem (TSP) asks for a minimum-cost Hamiltonian cycle
over a finite set of locations. The statement is short, but the feasible region
grows factorially with the number of cities, which has made the TSP a canonical
NP-hard reference in operations research and computer science. Its continued
relevance is equally practical: the same sequencing structure governs vehicle
and delivery routing in logistics, job and tool sequencing in manufacturing,
inspection and survey planning, communication network design, and motion
planning for robotic systems, including path generation for parallel
manipulators \cite{mottola2025applications}. Recent applications include
coordinated UAV routing and obstacle-aware mission planning
\cite{athira2024acodtsp,debnath2024uav}, while extensions such as the min--max
multiple TSP distribute workload among several agents \cite{perez2026partition}.
In these settings the operational question is rarely whether an optimal tour
exists; it is how close to one a method can come within a computational budget
that is fixed in advance. Exact algorithms remain indispensable when an
optimality certificate is required, and polyhedral analysis supports their
extension to related network design variants \cite{labbe2004ringstar}, but their
worst-case cost restricts routine use as instance size or operational complexity
increases, while alternative computational paradigms such as quantum search
formulations remain exploratory \cite{sato2025twostep}. Metaheuristics accept the
loss of the certificate in exchange for high-quality solutions obtained under
finite budgets by combining diversification and intensification
\cite{talbi2009metaheuristics}. TSPLIB supplies the standardized instances and
distance conventions that make progress in this direction reproducibly
comparable \cite{reinelt1991}.

The breadth of recent TSP research, including tree-seed search \cite{cinar2020},
social-spider optimization \cite{bay2021}, discrete Jaya \cite{gunduz2021},
fish-migration optimization \cite{chu2020dfmo}, probe-based search
\cite{rahman2024probe}, and machine-learning-assisted solvers
\cite{alanzi2025tsp,sui2025survey}, shows that generating feasible tours is no
longer the binding constraint on performance. What limits a modern method is
coordination. A competitive design must simultaneously preserve permutation
feasibility under every operator it applies, convert the information accumulated
by the population into movement that is genuinely informative, retain enough
diversity to postpone premature agreement, provide local improvement strong
enough to repair edge structure that population operators repair only by chance,
keep the cost of that improvement bounded relative to the available budget, and
behave stably as the number of cities grows. These requirements are not
independent, and they are frequently antagonistic: strengthening intensification
consumes the budget that exploration needs, whereas preserving diversity defers
the exploitation on which final tour quality depends. The question addressed in
this paper is therefore not which operators to use, but how much effort each
level of the search should receive, how often, and on which solutions.

Particle Swarm Optimization (PSO) is an attractive starting point because it
makes population-level learning explicit: each particle combines its own search
history with information shared by the swarm, and the relative weight of the two
sources is a control rather than an emergent property
\cite{kennedy1995,eb1995,gad2022pso}. PSO was, however, defined over continuous
vector spaces. A tour is a permutation, so vector addition and velocity updates
have no literal counterpart; movement must be redefined through operators that
map permutations to permutations while preserving the adjacency information that
gives a tour its quality. Discrete formulations that treat this as a notational
detail tend to fail in characteristic ways: they generate moves that must be
repaired to remain feasible, they destroy high-quality edges faster than the
swarm can rediscover them, or they contract the population onto a single
attractor within a few iterations. Work on adjustable controls and
heterogeneous or self-adaptive parameters confirms both the promise of
discrete PSO and the sensitivity of its exploration dynamics to these choices
\cite{strak2018adjustability,strak2019selfadaptive}. It
also indicates that mutation-driven movement alone rarely delivers competitive
tour quality: the swarm layer decides where to look, but it does not by itself
extract the quality available in a promising region.

Local search supplies that missing capability. Swap, shift, insertion,
inversion, symmetry and 2-opt moves act directly on tour structure, and 2-opt in
particular removes the crossing edges that permutation mutation eliminates only
incidentally; this effectiveness explains its recurrence across DPSO, ACO, Jaya,
fish-migration and hybrid schemes
\cite{bas2023symmetry,gunduz2021,chu2020dfmo,uddin2023improvement}. The
difficulty is that the capability is not free and its cost scales unfavorably.
An unrestricted 2-opt pass inspects a quadratic number of city pairs, so
applying it to every particle at every iteration transfers the run from the
swarm to the neighborhood scan and drives the population toward a set of
mutually similar 2-optimal tours; applying it too rarely has the opposite
effect, leaving a diverse population and an unpolished final tour. A hybrid that
resolves this tension must answer four separable questions: which solutions
receive intensification, which moves are admissible, how often refinement is
triggered, and how much of the strongest refinement is reserved for the
incumbent. It must answer them through controls that can be reported,
tuned and audited.

ADPSO-ERLS answers these four questions with three coordinated mechanisms and
one accounting rule, which together constitute the algorithmic novelty of this
work:

\begin{itemize}
\item \textbf{Memory-driven permutation movement with heterogeneous
seeding.} Particles move by imitating either their personal best or the
incumbent and applying a single transposition, so feasibility never has to be
repaired; a tuned fraction of the swarm is seeded by nearest-neighbor
construction while the remainder stays uniformly random, which raises initial
quality without homogenizing the population.

\item \textbf{Elitist candidate-restricted local search (ERLS).}
During evolution, 2-opt is confined to the $K$-nearest candidate lists, applied
only to an elite subset of the swarm, only every $\ell$ iterations, and only for
at most $T_2$ passes. Which solutions are refined, which exchanges are
admissible, how often refinement occurs, and how much work each call may perform
are thus four separate inputs rather than a fixed schedule.

\item \textbf{Staged final refinement with controlled
perturbation.} The strongest effort is concentrated on the incumbent alone
through candidate-restricted 2-opt, an optional full 2-opt phase, and
double-bridge perturbations followed by candidate-restricted repair, each
committed only on strict improvement, so that
routine exploitation and expensive end-of-run intensification are separately
budgeted instead of being merged into a single operator applied at constant
intensity.

\item \textbf{Recorded evaluation accounting and documented tuning.}
The swarm, local-search, and final-refinement stages use one candidate-assessment
counter. Elapsed time is reported alongside quality because equal assessment
counts do not imply equal arithmetic work for different operators.
\end{itemize}

The review below examines how far existing work has gone in that direction, and
Section~\ref{sec:contributions} states the resulting contributions of the study
in full.

The relationship with classical PSO should be stated precisely. ADPSO-ERLS has
no velocity vector, inertia weight, or acceleration coefficients. A particle
selects a personal or global memory, generates up to $s$ independent swap
mutants from that memory, retains the lowest-cost mutant, and replaces its
current position unconditionally; only personal and global memories are updated
elitistically. The method is therefore PSO-inspired through its population and
personal/global memory structure, not through a discrete analogue of the
classical velocity equation.

\subsection{Related works}\label{sec:related}
The literature bearing on this question can be read as five methodological
lines, each answering part of the coordination problem and each leaving a
different part open.

\subsubsection{Discrete PSO and permutation learning}

This line addresses the first question: how a swarm can move between feasible
permutations while still transmitting the information held in its memories. The
mechanism is the replacement of continuous velocities by combinatorial
operators. Early discrete formulations demonstrate the generality of the idea
beyond routing \cite{izakian2010}. The contribution of this group is
a feasible and interpretable notion of movement; its limitation is that the
intensity of that movement is itself a design variable which such formulations
largely fix in advance. A short operator sequence preserves valuable edges but
advances slowly, whereas an aggressive transformation erases precisely the
structural information the swarm was intended to share.

Adaptive control responds to this sensitivity. Strak et al. studied
adjustability, pheromone memory, heterogeneous parameter assignment and
self-adaptation for dynamic TSPs
\cite{strak2018adjustability,strak2019selfadaptive}, showing that diversity-aware
control and the reuse of search information after a change both improve
robustness. Their target, however, is adaptation to a moving optimum rather than
the allocation of intensification effort on a static instance, and their
pheromone-based memory differs from the staged refinement considered here.
Discrete PSO has also been redesigned for scheduling and channel-allocation
problems \cite{wang2023dynamic,qin2024overlapping}, which confirms the
representational flexibility of the paradigm while underlining that operators
must be matched to the constraints of each domain. Across this group, local
improvement is either absent or attached as an auxiliary step, and its budget is
not treated as a control of the method.

\subsubsection{Local-search-enhanced population methods}

Methods in this group address the complementary question: how to extract quality
once a promising region has been reached. The usual mechanism is a battery of
permutation neighborhoods applied inside the population loop. The DPSO family of
Ba\c{s} and Y{\i}ld{\i}zdan combines swap, shift, symmetry and 2-opt moves, with
DPSO6 assembling all four \cite{bas2023symmetry}, while improved 2-opt variants
concentrate on making a single neighborhood more
effective \cite{uddin2023improvement}. Population methods outside PSO reach the
same design point by other routes: DJAYA transforms solutions toward the best and
away from the worst population member \cite{gunduz2021}, DFMO alternates global
and local modes with swap, inversion and insertion \cite{chu2020dfmo}, and
two-stage probe search couples filtered route construction with local
augmentation and a self-escape mechanism \cite{rahman2024probe}. Collectively,
these studies establish that intensification is what separates competitive from
mediocre tour quality. What they leave open is its accounting: when
neighborhoods are broad, applied uniformly across the population, or repeated at
every iteration, the effort they consume grows with the number of cities and is
seldom reported as a tunable quantity, so the trade-off between quality and cost
remains implicit in the implementation.

\subsubsection{Large-step perturbation and staged refinement}

A third line treats escape from local optima as a concern separate from routine
improvement. Large-step Markov chains established the principle that a strong,
structure-preserving perturbation followed by re-optimization explores the space
of local optima far more effectively than restarting or than weak random noise
\cite{martin1991large}, with the double-bridge move as its canonical instrument.
In a hybrid architecture the value of the principle depends on where the
perturbation is placed: small restricted moves suit recurring exploitation
inside the population loop, whereas strong kicks are expensive and are most
useful when reserved for stagnation or for a final phase applied to a single
solution. Architectures that apply one local-search operator at constant
intensity throughout the run cannot express this distinction, and therefore
either pay the cost of strong intensification without concentrating its benefit,
or forgo the escape mechanism altogether. This separation of roles is what the
staged design of ADPSO-ERLS makes explicit.

\subsubsection{Reinforcement learning and data-informed search}

A fourth line replaces hand-set control decisions with learned policies or
predictions. DPSO-Q couples Q-learning with a discrete swarm so that search
behavior adapts online \cite{kappagantula2025dpsoq}, while QL-ACO updates
Q-values alongside pheromone trails, applies 2-opt, and isolates the
reinforcement-learning contribution through an ablation study
\cite{mohammed2026qlaco}. Both are forms of online learning and should be
distinguished from the offline tuning used by ADPSO-ERLS, which selects an
instance configuration before the final independent runs and adapts nothing
during a run. The broader learning-based literature is wider still: surveys
separate constructive, improvement, direct-hybrid and emerging
language-model-based solvers, and identify generalization, interpretability and
attainable solution quality as the persistent open questions
\cite{alanzi2025tsp,sui2025survey}. GTG-ACO learns heuristic and initial
pheromone matrices with a graph transformer \cite{abir2025gtgaco}, and ACO-ML
injects edge-classification probabilities into initialization, clustering and
transition decisions \cite{stodola2026acoml}. Such methods can transfer
structural priors across instances, an ability the present work does not claim;
in exchange they introduce a training pipeline and a learned component whose
effort profile is harder to audit than an explicitly reported budget.

\subsubsection{Ant-colony and hybrid schemes}

ACO approaches the same coordination problem from a constructive angle: edge
desirability is learned through pheromone reinforcement and tours are built
probabilistically, which differs structurally from the personal- and
global-best permutation learning of DPSO. Parameter optimization improves
conventional ACO \cite{wang2021ant}, and adaptive node clustering combines
candidate groups, entropy-controlled evaporation and diversity-based termination
to counter stagnation and parameter sensitivity \cite{stodola2022clustering}.
AACO-LST adapts transition and pheromone rules and applies 2-opt to selected
tours for large-scale instances \cite{tang2023adaptive}, which is the same
selective-intensification principle exploited here, realized within a
constructive rather than a population-movement framework. Hybrid variants
illustrate further divisions of labor: bi-heuristic schemes use PSO or
flower-pollination mechanisms to supervise ACO parameters while 2-opt performs
local search \cite{rokbani2021biheuristic}, a scheme-library ACO reuses stored
information and 2-opt in dynamic environments \cite{wang2022scheme}, and AAACO
combines ACO seeding with artificial-algae refinement and diversity-preserving
permutation operators \cite{altiok2026aaaco}. Related routing variants extend
the same machinery to UAV workload constraints \cite{athira2024acodtsp} and to
balanced min--max multiple-TSP tours \cite{perez2026partition}. These results
confirm that selective refinement and adaptive control are effective, while also
showing that their benefit depends on how information transfer and neighborhood
effort are coordinated. This coordination is usually embedded in the construction
rule rather than exposed as an independent budget.

\subsubsection{Research gap}

Each ingredient relevant to ADPSO-ERLS is therefore already available in the
literature: discrete swarm learning, heterogeneous or learned control, candidate
restriction, 2-opt, large-step perturbation, and staged or hybrid improvement.
No novelty is claimed for any ingredient in isolation. The narrower design
question examined here is whether four decisions about in-run intensification
can be exposed independently: which solutions are refined, which moves are
admissible, how often refinement is invoked, and how many passes a call may
perform. Table~\ref{tab:positioning} records only those decisions. This is a
narrative comparison of the cited methods, not a systematic review; it therefore
supports the positioning of the present study but not an exhaustive priority
claim.

\begin{table}[H]
\centering\small
\caption{Four independently exposed in-run intensification decisions in the
cited methods. ``Not exposed'' means that the publication does not present the
decision as an independent control; it does not imply that the method lacks all
related search behavior.}
\label{tab:positioning}
\begin{tabular}{llcccc}
\toprule
Method & Reference & Selection & Restriction & Frequency & Pass budget\\
\midrule
SA-DPSO         & \cite{strak2019selfadaptive} & Not exposed & Not exposed & Not exposed & Not exposed\\
DPSO6           & \cite{bas2023symmetry}       & Not exposed & Yes & Not exposed & Not exposed\\
DPSO-Q          & \cite{kappagantula2025dpsoq} & Not exposed & Not exposed & Not exposed & Not exposed\\
AACO-LST        & \cite{tang2023adaptive}      & Yes & Yes & Not exposed & Not exposed\\
AAACO           & \cite{altiok2026aaaco}       & Not exposed & Not exposed & Not exposed & Not exposed\\
AACO-NC         & \cite{stodola2022clustering} & Yes & Yes & Not exposed & Not exposed\\
\midrule
\textbf{ADPSO-ERLS} & this work & Yes & Yes & Yes & Yes\\
\bottomrule
\end{tabular}
\end{table}

A second gap concerns evaluation. Studies frequently report a single quality
statistic or a narrow size range, which leaves unanswered whether an advantage
in best-case quality survives at the level of typical runs, whether it is stable
across repetitions, whether it is supported by paired and global statistical
evidence, at what computational cost it is obtained, and whether the method
remains usable at larger sizes. Reporting solution quality, repeatability,
statistical evidence, computational limitations and scalability jointly is what
makes it possible to judge when additional refinement effort is worthwhile, and
this study is organized to provide that joint account.

These gaps define the paper's positioning with respect to the reviewed
literature. ADPSO-ERLS is not proposed as another collection of permutation
operators, nor as a replacement for specialized TSP solvers. It is a
budget-explicit DPSO framework in which recurring restricted improvement,
incumbent-focused final refinement and perturbation are separated into auditable
levels of effort. This separation permits the quality gained by hybridization
to be discussed together with its computational cost, which is obscured when
local search is embedded in a fixed or undocumented schedule.

\subsection{Our contributions}\label{sec:contributions}
ADPSO-ERLS is designed around these two gaps. Its organizing idea is to
stratify the search by level of effort rather than to add a further operator.
Exploration is carried out by a discrete swarm whose particles move by
memory-driven permutation transformations, so that feasibility is preserved by
construction. Intensification during the run is deliberately restricted:
candidate lists limit which exchanges are admissible, an elite subset limits
which solutions are refined, an interval limits how often refinement occurs, and
a pass budget limits how much work each call may perform. The strongest
intensification is deferred to a staged final refinement applied to the
incumbent alone, in which candidate-restricted 2-opt, an optional full 2-opt
phase, and LKH-like perturbations are executed in sequence and committed only on
strict improvement. The term LKH-like refinement denotes double-bridge kicks
alternated with cyclic candidate-restricted repair, not an implementation of the
variable-depth exchanges of the Lin--Kernighan--Helsgaun heuristic. Each level
therefore operates at the cost that its expected return justifies. The claim of
this paper is not that discrete PSO, 2-opt, candidate lists, double-bridge kicks
or automated tuning are new, but that coordinating them through explicit,
separately tunable budgets is a design decision worth isolating, reporting and
evaluating.

The empirical study assesses this design over 50 independent runs per instance
on five symmetric TSPLIB instances against five comparison implementations, six
methods in total. ADPSO-ERLS attains the lowest best cost and the lowest mean
cost on all five instances, with Gap values between 1.93\% and 4.17\%, RE values
between 3.04\% and 5.50\%, and APD values between 0.95\% and 1.28\%, so that the
advantage holds for a typical run rather than for a single fortunate execution.
The dispersion results reinforce this conclusion: the reported standard
deviation is the smallest for d493 and the second smallest for d657, rat783,
pr1002 and u1060. Thus, ADPSO-ERLS combines the lowest mean cost with a first or
second position in run-to-run variability on every comparison instance.
The Friedman test assigns it mean rank 1.0, the Iman--Davenport correction
confirms the global difference, and all 25 paired Wilcoxon comparisons
favor it with effect sizes classified as large. A separate study establishes feasible operation from
1,400 to 16,862 cities, with Gap values below 6.7\% under the recorded
configurations. All six methods are implemented in Rust and share the strict stopping value
$B=\texttt{max\_fes}=100000$, holding the evaluation budget and programming language common across methods and
making runtime a relevant empirical dimension.

Specifically, this work provides the following contributions.
\begin{itemize}
\item \textbf{A permutation-safe hybrid swarm architecture.} ADPSO-ERLS retains
personal- and global-best learning while moving through tour-preserving discrete
operators, so that every solution produced by the swarm is a valid tour without
repair. This resolves the representational mismatch that limits direct transfers
of PSO to the TSP, and it matters because feasibility maintained by construction
allows the remaining design effort to be spent on where and how strongly to
search. The architecture is documented in Section~\ref{sec:architecture} and
stated compactly in Algorithm~\ref{alg:main}.

\item \textbf{Restricted, selective, scheduled and budgeted in-run
intensification.} Candidate lists bound the neighborhood inspected by local
search, an elite subset determines which solutions are refined, and the
local-search interval and 2-opt pass budget determine how often and how strongly
exploitation acts during swarm evolution. This addresses the scaling tension
identified above, since it prevents neighborhood scanning from absorbing the run
while retaining the corrective power of edge exchange. All four quantities are
independent inputs listed in Table~\ref{tab:parameters} and tuned rather than
fixed.

\item \textbf{A staged final refinement with controlled perturbation.} Subject
to the remaining global evaluation budget, the incumbent receives candidate-restricted 2-opt,
an optional full 2-opt phase, and LKH-like perturbations, each committed only on
strict improvement, which makes the cascade monotone by construction. This
separates routine exploitation from expensive end-of-run effort and concentrates
the latter where a single high-quality solution can absorb it. Its cost is
charged explicitly to the same strict FES budget as initialization and swarm
evolution.

\item \textbf{Reproducible hyperparameter control.} The principal exploration
and refinement variables are exposed in the documented Optuna search space of
Table~\ref{tab:searchspace} and the selected values are reported per instance in
Table~\ref{tab:optuna_configurations}. This replaces undocumented implementation
constants by an auditable selection procedure, so that the balance between
quality and computational cost can be inspected and reproduced; the resulting
variation across instances is reported as an observed association, not as a
causal effect of instance size.

\item \textbf{A multi-dimensional and statistically supported evaluation.} Five
TSPLIB instances are evaluated over 50 independent runs against five comparison
methods, reporting best-case accuracy (Gap), expected accuracy (RE), internal
repeatability (APD), standard deviation, descriptive statistics and confidence intervals together,
and combining a six-method Friedman analysis with 25 paired Wilcoxon
comparisons. This prevents the single-best-tour reporting that cannot
distinguish a consistently strong method from a fortunate one, and it separates
statistical significance from the magnitude of the cost differences.

\item \textbf{Scalability evidence accompanied by explicit limitations.}
ADPSO-ERLS is applied to seven further instances of 1,400 to 16,862 cities,
establishing practical feasibility at these sizes while identifying quadratic
distance storage and intensive local-search passes as the operative constraints.
The measurements are presented as evidence of feasibility under the common
budget $B=100000$, not as an asymptotic scaling law, and runtime is
interpreted descriptively for the reasons stated above.
\end{itemize}
Taken together, these contributions position ADPSO-ERLS as a reproducible
framework for studying how discrete swarm exploration and progressively stronger
local refinement can be coordinated under explicit computational controls.

The remainder of the paper is organized as follows.
Section~\ref{sec:preliminaries} fixes the notation and the distance conventions
used throughout. Section~\ref{sec:algorithm} describes the architecture: its
control parameters, its operators, the three operational phases and the
evaluation-accounting rule, and states the method as a single procedure in
Algorithm~\ref{alg:main}. Section~\ref{sec:experimental-results} reports the
Optuna study that selects the instance configurations, the comparison against
five baselines over 50 independent runs on five TSPLIB instances, the Friedman
and Wilcoxon analyses, and the scalability experiments up to 16,862 cities.
Section~\ref{sec:conclusions} summarizes what the evidence supports, delimits
what it does not, and derives a prioritized agenda from the remaining
limitations.

\section{Preliminaries}\label{sec:preliminaries}
\subsection{Mathematical formulation and problem definition}

\begin{definition}[Traveling Salesman Problem (TSP)]\label{def:tsp}
	Given a complete weighted graph $G = (V, E, d)$ where:
	\begin{itemize}
		\item $V = \{1, 2, \dots, m\}$ is the vertex set representing $m$ cities.
		\item $E$ is the set of all edges between pairs of vertices (complete graph).
		\item $d: E \to \mathbb{R}^+$ is a weight function assigning distances or costs between cities, specifically,
		\begin{equation}\label{eq:euclidean_distance}
			d_{ij} = \|\mathbf{c}_i - \mathbf{c}_j\|_2 = \sqrt{(x_i - x_j)^2 + (y_i - y_j)^2},
		\end{equation}
		represents Euclidean distances between cities $i$ and $j$. Moreover, $\mathbf{c}_i = (x_i, y_i) \in \mathbb{R}^2$ are the Cartesian coordinates of city $i$.
	\end{itemize}
	The objective is to find a Hamiltonian cycle $\pi$ minimizing the total tour length:
	\begin{equation*}
		\min_{\pi \in \mathcal{S}_m} f(\pi),
	\end{equation*}
	where
	\begin{equation}\label{eq:tsp_objective}
		f(\pi) = \sum_{i=1}^{m-1} d_{\pi(i),\pi(i+1)} + d_{\pi(m),\pi(1)},
	\end{equation}
	and $\pi: \{1,\dots,m\} \to \{1,\dots,m\}$ is a bijective permutation function representing the order of city visits. Here, $\mathcal{S}_m$ denotes the set of all permutations of $\{1,2,\dots,m\}$.
\end{definition}

The TSP constitutes an NP-hard combinatorial optimization problem characterized by a factorial growth in solution space complexity, with $|\mathcal{S}_m| = m!$ possible tours for $m$ cities. The Euclidean variant considered here satisfies the triangle inequality $d_{ij} \leq d_{ik} + d_{kj}$ for all $i,j,k \in V$, which enables certain algorithmic optimizations but does not reduce the fundamental computational complexity.

\subsubsection{Distance matrix construction}

In practical implementations, TSP instances are typically provided in standardized formats such as TSPLIB \cite{reinelt1991}, containing the Cartesian coordinates $\mathbf{c}_i = (x_i, y_i)$ for each city $i \in \{1,\dots,m\}$. The algorithm initiates with a critical preprocessing phase that computes the symmetric distance matrix $\mathbf{D} \in \mathbb{R}^{m \times m}$,
\begin{equation}\label{matrixD}
	\mathbf{D} =\begin{pmatrix}
		d_{11} & \cdots & d_{1m}\\
		\vdots & \ddots & \vdots\\
		d_{m1} & \cdots & d_{mm}
	\end{pmatrix},
\end{equation}
where each entry $d_{ij}$ is defined according to \eqref{eq:euclidean_distance}.

This matrix exhibits the following mathematical properties that are fundamental to the Euclidean TSP:
\begin{itemize}
	\item Symmetry: $d_{ij} = d_{ji}$ for all $i,j \in V$, reflecting the undirected nature of the distance metric.
	\item Non-negativity: $d_{ij} \geq 0$ with $d_{ii} = 0$, ensuring all distances are positive and self-distances are zero.
	\item Triangle inequality: $d_{ij} \leq d_{ik} + d_{kj}$ for all $i,j,k \in V$, a geometric property that enables algorithmic optimizations and informs local search operators.
\end{itemize}

\section{Algorithm description}\label{sec:algorithm}

Throughout this section, the problem is the one fixed in
Definition~\ref{def:tsp}, and tour costs are computed under the
\texttt{EUC\_2D} convention of Eq.~\eqref{eq:euclidean_distance}. The
description below makes the computational budget part of the algorithm rather
than an implementation detail. In particular, initialization, evolutionary
search, in-run local search, and final refinement all consume units from the
same global counter, and no candidate may be assessed once its active budget
limit has been reached.

\subsection{Algorithm overview}

ADPSO-ERLS couples a discrete particle-swarm layer, which determines where to
search, with a local-search layer, which determines how much quality to extract
from selected regions. In the continuous PSO of Kennedy and Eberhart
\cite{kennedy1995,eb1995}, a particle moves by means of a velocity combining
inertia and attraction toward personal and global memories. A TSP tour is a
permutation rather than a point in a Euclidean vector space, so ADPSO-ERLS
replaces the velocity equation by a permutation-safe, memory-driven
transformation. At each update, particle~$i$ selects its personal best
$\mathit{pbest}_i$ with probability~$\gamma$ and the global best~$\pi_g$ with
probability~$1-\gamma$, and then applies a random swap to a copy of the selected
memory. Among at most~$s$ independently generated swap candidates, the
shortest one becomes the particle's new position.

Two additional mechanisms control the balance between exploration and
intensification.
\begin{itemize}
	\item[\checkmark] \textbf{Heterogeneous initialization.} An elite fraction
	$\alpha$ of the swarm is initialized by nearest-neighbor construction from
	random starting cities, whereas the remaining particles are random
	permutations. This provides high-quality initial tours without making the
	whole population structurally homogeneous.
	
	\item[\checkmark] \textbf{Budgeted two-level local search.}
	Candidate-restricted 2-opt is applied periodically to the elite subset during
	evolution. A stronger incumbent-only cascade is then executed within a
	budget reserved in advance. The cascade combines candidate-restricted 2-opt,
	optional full 2-opt, and cyclic candidate repair after double-bridge
	perturbations. Both levels use the same evaluation counter as the swarm.
\end{itemize}

The central scheduling rule is
\begin{equation}\label{eq:budget_partition}
	B_{\mathrm{evo}}=\left\lfloor\eta B\right\rfloor,
	\quad
	B_{\mathrm{fin}}=B-B_{\mathrm{evo}},
	\quad 0<\eta<1,
\end{equation}
where $B$ is the global evaluation limit and $\eta$ is the fraction assigned to
initialization and evolutionary search. Throughout this work $\eta=0.7$, so that
$B_{\mathrm{evo}}=70000$ and $B_{\mathrm{fin}}=30000$. The value of~$\eta$ is fixed
before hyperparameter tuning and is held constant across all instances, runs, and
ablation variants. Consequently, the final cascade always has an explicit
nonzero reserve and is never appended as uncounted post-processing.

\subsubsection{Algorithm parameters and outputs}

\paragraph{Inputs.}
A run is determined by the coordinates
$\mathbf c_1,\ldots,\mathbf c_m$, the parameters in
Table~\ref{tab:parameters}, and one integer seed. The seed initializes the
single pseudorandom generator used for nearest-neighbor starting cities, random
permutations, memory choices, swap positions, and double-bridge cut points.
The feasibility condition $P\leq B_{\mathrm{evo}}$ is imposed so that the
initial population can always be evaluated completely.

\paragraph{Control parameters.}
The elite cardinality is
$E=\max(1,\lceil\alpha P\rceil)$, and the effective candidate-list length is
$\min(K,m-1)$. The tolerance $\varepsilon=10^{-10}$ in
Eq.~\eqref{eq:2opt_gain} is a fixed numerical safeguard rather than a tuned
parameter. The global budget~$B$, the evolutionary share~$\eta$, and the
equal-third final-stage deadlines defined below are fixed by the experimental
protocol and are not optimized per instance.

\begin{table}[H]
	\centering\small
	\setlength{\tabcolsep}{6pt}
	\caption{Parameters of ADPSO-ERLS. Symbols are grouped according to the
		component they control. The strict global budget $B$ and the evolutionary
		share $\eta$ are protocol constants rather than tuned hyperparameters.}
	\label{tab:parameters}
	\begin{tabular}{@{}ll>{\raggedright\arraybackslash}p{8.1cm}@{}}
		\toprule
		\textbf{Symbol} & \textbf{Type} & \textbf{Description} \\
		\midrule
		\multicolumn{3}{@{}l}{\emph{Swarm}}\\
		$P$ & Integer & Number of particles \\
		$\alpha$ & Real in $[0,1]$ & Elite fraction, with
		$E=\max(1,\lceil\alpha P\rceil)$ \\
		$\gamma$ & Real in $[0,1]$ & Probability of selecting the personal best
		rather than the global best as mutation source \\
		$s$ & Integer & Maximum number of swap candidates assessed per particle
		update \\
		\midrule
		\multicolumn{3}{@{}l}{\emph{In-run local search}}\\
		$K$ & Integer & Number of nearest neighbors retained per city \\
		$\ell$ & Integer & Number of iterations between successive elite local-search
		calls \\
		$T_2$ & Integer & Maximum number of candidate-restricted 2-opt passes per
		elite refinement \\
		\midrule
		\multicolumn{3}{@{}l}{\emph{Final incumbent refinement}}\\
		$T_f$ & Integer & Maximum number of candidate-restricted 2-opt passes \\
		$T_F$ & Integer & Maximum number of full 2-opt passes; $T_F=0$ disables
		this stage \\
		$\kappa$ & Integer & Maximum number of double-bridge kicks \\
		$\mu$ & Integer & Maximum number of accepted improving moves per cyclic
		candidate repair \\
		\midrule
		\multicolumn{3}{@{}l}{\emph{Budget and reproducibility}}\\
		$B$ & Integer & Strict global limit on candidate-solution evaluations;
		$B=100000$ in the experiments \\
		$\eta$ & Real in $(0,1)$ & Fixed fraction of $B$ assigned to
		initialization and evolution; $\eta=0.7$ in the experiments \\
		--- & Integer & Base random seed \\
		\bottomrule
	\end{tabular}
\end{table}

\paragraph{Outputs.}
A run returns the incumbent tour $\pi_g\in\mathcal S_m$ and its stored cost
$f(\pi_g)$. Personal and global memories are replaced only after a strict
improvement. Likewise, each stage of the final cascade is committed only if it
improves the incumbent. The incumbent cost is therefore non-increasing, and the
returned tour is the shortest candidate accepted during the run. For
reproducibility, the implementation also records the seed, final counter, and
counter values at the end of initialization, evolution, and each final stage.

\subsubsection{Evaluation unit and strict accounting}\label{sec:fes-accounting}

A function-evaluation unit (FES) is charged whenever the algorithm determines
the objective value of a new candidate tour for the purpose of accepting,
rejecting, ranking, or storing it. This definition is independent of the
calculation used to obtain that value:
\begin{itemize}
	\item evaluating a complete tour through Eq.~\eqref{eq:tsp_objective} costs
	one FES;
	\item assessing a swap candidate costs one FES;
	\item assessing a 2-opt candidate costs one FES even when its cost is
	obtained in $O(1)$ time from an edge delta;
	\item assessing a double-bridge candidate costs one FES, and every subsequent
	2-opt candidate considered during its repair costs one additional FES.
\end{itemize}
Distance-matrix construction, candidate-list construction, random-number
generation, permutation copying, sorting particle indices, and updating a
stored cost after an already charged assessment do not create or assess a new
candidate and therefore do not increment the counter.

All search routines use the following atomic rule, where $L$ is the active
phase deadline:
\begin{equation}\label{eq:fes_invariant}
	\textsc{Assess}(z,L):\quad
	\begin{cases}
		\text{return }\textsc{BudgetExhausted}\text{ without assessing }z,
		& q\geq L,\\
		\text{determine }f(z)\text{; set }q\gets q+1\text{; return }f(z),
		& q<L.
	\end{cases}
\end{equation}
The determination of $f(z)$ in Eq.~\eqref{eq:fes_invariant} may be a complete
summation or an exact incremental calculation. The increment occurs immediately
after the candidate cost is determined and before any acceptance or memory
update. Thus, no routine can overshoot its deadline, and no assessed candidate
is hidden from the common counter.

This candidate-based definition is essential for a hybrid method: counting only
complete recomputations would assign zero FES to most 2-opt work merely because
the implementation uses exact edge deltas. Wall-clock time is nevertheless
reported separately because one FES can have different arithmetic cost in a
full evaluation and in an incremental neighborhood assessment.

\subsubsection{Fundamental algorithmic operators}\label{sec:operators}

The algorithm is assembled from standard operators; the contribution lies in
their budgeted combination and scheduling.

\begin{itemize}
	\item[\checkmark] \textbf{Distance matrix and candidate lists.}
	All pairwise distances are computed once and stored in the symmetric matrix
	$\mathbf D$ of Eq.~\eqref{matrixD}. For each city~$v$, the ordered list
	$\mathcal N_K(v)$ contains its $\min(K,m-1)$ nearest neighbors. The restricted
	2-opt variants use these lists, whereas full 2-opt enumerates unrestricted
	edge pairs.
	
	\item[\checkmark] \textbf{Nearest-neighbor construction.}
	Starting from a randomly selected city, the construction repeatedly appends
	the nearest unvisited city. It first scans the candidate list and falls back
	to all unvisited cities when necessary. Each completed initial tour is
	assessed once.
	
	\item[\checkmark] \textbf{Swap mutation.}
	Given a source tour~$\sigma$, two distinct positions are drawn uniformly and
	their cities are exchanged. Each resulting permutation is a candidate and is
	charged one FES before it can participate in selection.
	
	\item[\checkmark] \textbf{Elite selection and memory update.}
	Particles are ranked by $f(\mathit{pbest}_i)$, and the best~$E$ indices form
	the elite set~$\mathcal E$. A candidate replaces $\mathit{pbest}_i$ only if it
	is strictly shorter; it replaces~$\pi_g$ only if it also improves the current
	global record.
	
	\item[\checkmark] \textbf{2-opt exchange.}
	For tour edges $(a,b)$ and $(c,e)$, a 2-opt move deletes both edges, reverses
	the intervening segment, and reconnects the tour using $(a,c)$ and $(b,e)$.
	Its exact cost variation is
	\begin{equation}\label{eq:2opt_gain}
		\Delta=d_{ac}+d_{be}-d_{ab}-d_{ce}.
	\end{equation}
	Before computing and testing a candidate move, the routine applies
	Eq.~\eqref{eq:fes_invariant}. If budget remains, the move is charged one FES
	and accepted under first improvement whenever $\Delta<-\varepsilon$.
	Candidate-restricted 2-opt obtains possible reconnections from
	$\mathcal N_K(\cdot)$; full 2-opt considers all admissible edge pairs; and
	cyclic candidate 2-opt treats the tour circularly and bounds a call by at most
	$\mu$ accepted improvements. Every attempted move, whether improving or not,
	counts once because each one is an assessed candidate.
	
	\item[\checkmark] \textbf{Double-bridge perturbation.}
	Three cut points divide the incumbent into four segments, which are
	reassembled in a different order to reach another search basin
	\cite{martin1991large}. The perturbed tour is assessed once and is then
	repaired by budgeted cyclic candidate 2-opt. Only a repaired tour that
	strictly improves the incumbent is committed.
\end{itemize}

\subsection{Detailed architecture and budget flow}\label{sec:architecture}

Let $q$ denote the single FES counter. Preprocessing does not modify~$q$ because
it produces problem data rather than candidate tours. The search then executes
the following three phases.

\medskip
\noindent\textbf{Phase 1: heterogeneous initialization.}
The algorithm creates $E$ nearest-neighbor tours and $P-E$ random tours. Each
completed tour is assessed through Eq.~\eqref{eq:fes_invariant} with deadline
$B_{\mathrm{evo}}$, its personal best is initialized to itself, and the shortest
initial tour becomes~$\pi_g$. Consequently, the phase terminates with
\begin{equation}\label{eq:q_initialization}
	q_{\mathrm{init}}=P\leq B_{\mathrm{evo}}.
\end{equation}

\medskip
\noindent\textbf{Phase 2: evolutionary search and selective intensification.}
The evolutionary loop is executed while $q<B_{\mathrm{evo}}$. Every~$\ell$
iterations, the elite set is selected and each elite personal best is refined
by at most~$T_2$ candidate-restricted 2-opt passes. The routine stops immediately
when it reaches a local optimum, exhausts its pass allowance, or reaches
$B_{\mathrm{evo}}$. The refined tour replaces the particle's current position,
and memory updates are performed only after strict improvement.

The memory-driven update then processes particles in a fixed index order. For
each particle, personal or global memory is selected according to~$\gamma$, and
up to~$s$ swap candidates are generated and assessed. The best assessed mutant
becomes the new position even if it is worse than the previous position; only
the personal and global memories are elitist. If the active deadline is reached
part-way through a particle update, the best candidates already assessed for
that particle remain valid, no further candidate is generated, and all
remaining particles are left unchanged. Therefore, Phase~2 terminates without
overshoot at
\begin{equation}\label{eq:q_evolution}
	q_{\mathrm{evo}}=B_{\mathrm{evo}},
\end{equation}
provided that at least one swap candidate can be generated, as is the case for
$m\geq2$ and $s\geq1$.

\medskip
\noindent\textbf{Phase 3: reserved final refinement.}
The final reserve is partitioned by three cumulative deadlines:
\begin{equation}\label{eq:final_deadlines}
	L_1=B_{\mathrm{evo}}+\left\lfloor\frac{B_{\mathrm{fin}}}{3}\right\rfloor,
	\quad
	L_2=B_{\mathrm{evo}}+\left\lfloor\frac{2B_{\mathrm{fin}}}{3}\right\rfloor,
	\quad
	L_3=B.
\end{equation}
Candidate-restricted 2-opt runs first and may assess candidates only while
$q<L_1$. Full 2-opt then runs, when $T_F>0$, only while $q<L_2$. Finally, the
cyclic candidate and double-bridge stage runs only while $q<L_3$. These are
cumulative rather than isolated quotas: if an earlier stage converges before
its deadline or is disabled, its unused evaluations are inherited by later
stages. Conversely, an earlier stage cannot consume the nominal minimum
opportunity reserved for its successors.

Each stage starts from the current incumbent and is interrupted cleanly if its
deadline is reached. An incomplete neighborhood scan does not invalidate
previously accepted improvements, but no unassessed move is accepted. A
double-bridge kick is not performed unless at least one FES remains to assess
the perturbed tour. If the budget expires during its repair, the best feasible
tour already assessed in that repair is compared with the incumbent. The
commit-on-improvement rule guarantees that final refinement cannot deteriorate
$\pi_g$.

Table~\ref{tab:fes_phase_trace} gives the counter invariant recorded by the
implementation. The tuple
$(q_{\mathrm{init}},q_{\mathrm{evo}},q_1,q_2,q_3)$ is written to the run log,
which permits the budget flow of every execution to be audited directly.

\begin{table}[H]
	\centering\small
	\caption{FES trace and hard counter bounds for one run. The actual values
		$q_1,q_2,q_3$ depend on convergence and the pass limits, whereas no value may
		exceed its corresponding cumulative deadline.}
	\label{tab:fes_phase_trace}
	\begin{tabular}{@{}lccc@{}}
		\toprule
		\textbf{Stage} & \textbf{Initial counter} & \textbf{Final counter} &
		\textbf{Hard bound} \\
		\midrule
		Initialization & $0$ & $q_{\mathrm{init}}=P$ & $P\leq B_{\mathrm{evo}}$ \\
		Evolution and elite 2-opt & $P$ & $q_{\mathrm{evo}}$ &
		$q_{\mathrm{evo}}=B_{\mathrm{evo}}$ \\
		Final candidate 2-opt & $q_{\mathrm{evo}}$ & $q_1$ & $q_1\leq L_1$ \\
		Final full 2-opt & $q_1$ & $q_2$ & $q_2\leq L_2$ \\
		Cyclic repair and kicks & $q_2$ & $q_3$ & $q_3\leq L_3=B$ \\
		\bottomrule
	\end{tabular}
\end{table}

The final counter may be smaller than~$B$ if all enabled final routines reach
their stopping conditions before exhausting the reserve. This does not violate
the protocol: $B$ is a common maximum, not a requirement to perform redundant
assessments after convergence. The invariant required for every run is
\begin{equation}\label{eq:strict_budget}
	P=q_{\mathrm{init}}\leq q_{\mathrm{evo}}=B_{\mathrm{evo}}\leq q_1\leq q_2\leq q_3\leq B,
\end{equation}
in agreement with the exact termination $q_{\mathrm{evo}}=B_{\mathrm{evo}}$
established in Eq.~\eqref{eq:q_evolution} and recorded in
Table~\ref{tab:fes_phase_trace}.

\subsubsection{Complexity and memory}

Precomputing~$\mathbf D$ requires $O(m^2)$ time and memory, and candidate lists
require $O(mK)$ additional storage. A complete tour evaluation and a tour copy
are $O(m)$. A swap-cost implementation may use either a complete evaluation or
an exact incremental calculation, but both consume one FES. A
candidate-restricted 2-opt pass inspects at most $O(mK)$ moves, whereas an
unrestricted pass inspects $O(m^2)$ moves in the worst case; each inspected
move consumes one FES until the relevant deadline is reached. Segment reversal
and double-bridge construction are $O(m)$. Peak memory is dominated by the
$m^2$ distance matrix.

\subsection{Algorithm pseudocode}

Algorithm~\ref{alg:main} states the complete budget flow. Every auxiliary
routine that assesses a candidate receives an active deadline and implements
Eq.~\eqref{eq:fes_invariant}; a routine returns immediately when that deadline
is reached. Thus, budget checks occur inside neighborhood and particle loops,
not only at the top of the evolutionary loop.

\begin{algorithm}[H]
	\caption{ADPSO-ERLS with strict global evaluation accounting}
	\label{alg:main}
	\begin{algorithmic}[1]
		\Require coordinates of $m$ cities; $P,\alpha,\gamma,s,K,\ell,T_2,T_f,T_F,
		\kappa,\mu$; global budget $B$; fixed evolutionary share $\eta$;
		seed
		\Ensure incumbent tour $\pi_g$ and stored cost $f(\pi_g)$
		\State $\mathbf D\gets$\Call{PrecomputeDistances}{coordinates}
		\State $\mathcal N_K\gets$\Call{BuildCandidateLists}{$\mathbf D,K$}
		\State $B_{\mathrm{evo}}\gets\lfloor\eta B\rfloor$;
		$B_{\mathrm{fin}}\gets B-B_{\mathrm{evo}}$
		\State $L_1\gets B_{\mathrm{evo}}+\lfloor B_{\mathrm{fin}}/3\rfloor$;
		$L_2\gets B_{\mathrm{evo}}+\lfloor2B_{\mathrm{fin}}/3\rfloor$;
		$L_3\gets B$
		\State $E\gets\max(1,\lceil\alpha P\rceil)$;
		$q\gets0$;
		$f(\pi_g)\gets+\infty$
		\State $(\Pi,\pi_g,q)\gets$\Call{InitPopulation}{$P,E,\mathbf D,
			\mathcal N_K,q,B_{\mathrm{evo}}$}
		\Comment{Phase 1; one FES per completed initial tour}
		\State $q_{\mathrm{init}}\gets q$;
		$t\gets0$
		\While{$q<B_{\mathrm{evo}}$}
		\If{$t\bmod\ell=0$}
		\State $\mathcal E\gets$\Call{SelectElite}{$\Pi,E$}
		\State \Call{EliteLocalSearch}{$\mathcal E,\Pi,\mathbf D,
			\mathcal N_K,T_2,q,B_{\mathrm{evo}}$}
		\Comment{one FES per attempted 2-opt move}
		\EndIf
		\If{$q<B_{\mathrm{evo}}$}
		\State \Call{UpdatePositions}{$\Pi,\mathbf D,\gamma,s,q,
			B_{\mathrm{evo}}$}
		\Comment{stop inside the particle loop at the deadline}
		\EndIf
		\State $t\gets t+1$
		\EndWhile
		\State $q_{\mathrm{evo}}\gets q$
		\State $(z,c,q)\gets$\Call{CandidateTwoOpt}{$\pi_g,\mathbf D,
			\mathcal N_K,T_f,q,L_1$}
		\If{$c<f(\pi_g)$}
		\State $\pi_g\gets z$;
		$f(\pi_g)\gets c$
		\EndIf
		\State $q_1\gets q$
		\If{$T_F>0 \land q<L_2$}
		\State $(z,c,q)\gets$\Call{FullTwoOpt}{$\pi_g,\mathbf D,T_F,q,L_2$}
		\If{$c<f(\pi_g)$}
		\State $\pi_g\gets z$;
		$f(\pi_g)\gets c$
		\EndIf
		\EndIf
		\State $q_2\gets q$
		\If{$q<L_3$}
		\State $(z,c,q)\gets$\Call{LKStyleRefinement}{$\pi_g,\mathbf D,
			\mathcal N_K,\kappa,\mu,q,L_3$}
		\If{$c<f(\pi_g)$}
		\State $\pi_g\gets z$;
		$f(\pi_g)\gets c$
		\EndIf
		\EndIf
		\State $q_3\gets q$
		\State \Call{LogBudgetTrace}{$q_{\mathrm{init}},q_{\mathrm{evo}},q_1,q_2,q_3$}
		\State \Return $\pi_g,f(\pi_g)$
	\end{algorithmic}
\end{algorithm}

The pseudocode deliberately leaves the standard neighborhood loops inside the
named auxiliary routines, but their accounting contract is explicit: before
the cost of any new swap, 2-opt, or perturbed candidate is determined, the
routine must verify $q<L$; after determining that cost, it must increment~$q$
exactly once. The returned trace therefore demonstrates whether the
implementation has respected Eq.~\eqref{eq:strict_budget} in every individual
run.

\section{Experimental results}
\label{sec:experimental-results}

This section evaluates ADPSO-ERLS from three complementary perspectives. First,
the experimental framework defines the benchmark instances, comparison methods,
evaluation measures, stopping conditions, and reproducibility limitations.
Second, the Optuna study documents the instance-level hyperparameter selection,
after which solution quality, robustness, runtime, and statistical evidence are
examined over 50 independent runs. Finally, the scalability experiments and
their best-tour visualizations assess the practical behavior of the method on
instances containing up to 16,862 cities.

\subsection{Experimental framework}
The evaluation uses the symmetric Euclidean TSPLIB instances d493, d657,
rat783, pr1002, and u1060, containing 493, 657, 783, 1002, and 1060 cities,
respectively. The reference optimal costs used by the experiment are 35002, 48912,
8806, 259045, and 224094.  We compare ADPSO-ERLS with SA-DPSO
\cite{strak2019selfadaptive}, AACO-LST \cite{tang2023adaptive}, DPSO6 \cite{bas2023symmetry}, DPSO-Q
\cite{kappagantula2025dpsoq}, and AAACO \cite{altiok2026aaaco}. All six implementations are in
Rust and were run on the HP Laptop 15-dy2xxx described in
Table~\ref{tab:computing_platform}.

\begin{table}[H]
\centering\small
\caption{Hardware and software platform used for the experiments.}
\label{tab:computing_platform}
\begin{tabular}{ll}
\toprule
Component & Specification \\
\midrule
Computer model & HP Laptop 15-dy2xxx \\
Processor & 11th Gen Intel Core i7-1165G7, 8 logical processors \\
Memory & 16.0 GiB \\
Storage capacity & 512.1 GB \\
Graphics & Intel Iris Xe Graphics (TGL GT2) \\
Operating system & Ubuntu 26.04 LTS, 64 bit \\
Linux kernel & 7.0.0-28-generic \\
Firmware & F.25 \\
Desktop and display server & GNOME 50 with Wayland \\
\bottomrule
\end{tabular}
\end{table}

We evaluate performance with three complementary metrics. For a given
method--instance pair, let $R=50$ be the number of independent runs and let
$f_r$ denote the tour cost returned in run $r$. We define
\begin{equation}\label{eq:sample-costs}
	f_{\min}=\min_{1\leq r\leq R} f_r,
	\qquad
	\bar f=\frac{1}{R}\sum_{r=1}^{R}f_r,
\end{equation}
and use $f^\star$ for the reference optimal cost of the corresponding TSPLIB
instance \cite{reinelt1991}. Gap measures the best observed solution quality
relative to $f^\star$ \cite{bas2023symmetry,bernardino2024family,
talbi2009metaheuristics,uddin2023improvement}:
\begin{equation}\label{eq:experimental-gap}
	\operatorname{Gap}(\%)=
	100\frac{f_{\min}-f^\star}{f^\star}.
\end{equation}
Relative Error (RE) measures the deviation of the mean outcome from the same
reference cost \cite{gunduz2021,wang2021ant}:
\begin{equation}\label{eq:experimental-re}
	\operatorname{RE}(\%)=
	100\frac{\bar f-f^\star}{f^\star}.
\end{equation}
Finally, Average Percentage Deviation (APD) measures the separation between
the mean outcome and the best solution returned by the method
\cite{tang2023adaptive,wang2021ant,wang2023dynamic}:
\begin{equation}\label{eq:experimental-apd}
	\operatorname{APD}(\%)=
	100\frac{\bar f-f_{\min}}{f_{\min}}.
\end{equation}
Lower values are better for all three metrics. Gap characterizes best-case
quality, RE characterizes typical quality, and APD characterizes run-to-run
consistency relative to the method's best observed solution.
We additionally report the standard deviation of the 50 tour costs as an
absolute measure of run-to-run dispersion. A smaller value indicates greater
repeatability, but it must be interpreted jointly with the mean because an
algorithm can be consistently far from the optimum.

\begin{table}[H]
\centering\small
\caption{Parameter list of the comparison algorithms.}
\label{tab:comparison_parameters}
\begin{tabularx}{\textwidth}{lX}
\toprule Algorithm & Parameters used in the experiments \\
\midrule
SA-DPSO \cite{strak2019selfadaptive} & $P=112$, neighborhood size $=30$, $B=100000$, $\rho=0.10$; heterogeneous $c_1,c_2,c_3,\omega$ distributions as specified by the source method.\\
AACO-LST \cite{tang2023adaptive} & $B=100000$, $m=1.5n$, candidate list $=20$, $A=2$, $B_a=3$, $\epsilon=0.1$, $\lambda=0.1$, $\rho_0=0.3$, $\omega=0.7$, $\gamma=0.8$, $s_0=30$, $Q=100$.\\
DPSO6 \cite{bas2023symmetry} & $P=100$, $B=100000$, $w=0.4$, $c_1=c_2=0.2$, candidate size $=20$, maximum velocity $=0.8$; swap, shift, symmetry, and final 2-opt.\\
DPSO-Q \cite{kappagantula2025dpsoq} & $P=150$, $B=100000$, learning rate $=0.7$, discount factor $=0.5$, exploration rate $=0.1$.\\
AAACO \cite{altiok2026aaaco} & $P=100$, $B=100000$, ACO seeding budget $=2.5\%$ of $B$, $\Delta=2$, $\epsilon=0.3$, $A_p=0.5$, $\alpha=1$, $\beta=5$, $\rho=0.5$, $Q=1$.\\
\bottomrule
\end{tabularx}
\end{table}

\subsection{Optuna-based hyperparameter selection}\label{sec:optuna}

The eleven instance-dependent hyperparameters
$\boldsymbol{\theta}=(\alpha,\gamma,P,K,s,\ell,T_2,T_f,T_F,\kappa,\mu)$
were selected offline through a separate Optuna study for each benchmark
instance. The global evaluation limit $B=100000$ and the evolutionary budget
share $\eta=0.7$ defined in Eq.~\eqref{eq:budget_partition} were fixed by the
experimental protocol and were excluded from the optimization domain. The
configuration retained for each instance was fixed before the 50 independent
runs reported in Section~\ref{sec:results}; the algorithm did not adapt these
hyperparameters during evaluation.

For a configuration $\boldsymbol{\theta}$ and seed $r$, let
$f_{\boldsymbol{\theta},r}$ denote the returned tour cost. Using the same
reference optimal cost $f^\star$ defined above, Optuna minimizes
\begin{equation}\label{eq:optuna-objective}
	J(\boldsymbol{\theta})
	=
	\frac{1}{3}
	\sum_{r\in\mathcal{S}_{\mathrm{tune}}}
	100\frac{f_{\boldsymbol{\theta},r}-f^\star}{f^\star}.
\end{equation}
Thus, the Python driver treats the compiled Rust solver as a black box and
evaluates each trial through three independent solver executions using
$\mathcal{S}_{\mathrm{tune}}=\{42,123,999\}$. Each instance-specific study
comprises 50 completed trials, corresponding to 150 solver executions. The
best completed trial according to Eq.~\eqref{eq:optuna-objective} was then
fixed for the final experiment. Final performance was assessed with the common
50-seed sequence used for all methods. All tuning executions followed the
strict global accounting rule in Eq.~\eqref{eq:strict_budget}: initialization,
evolutionary updates, local search, and final refinement drew from the same
limit $B$.

Table~\ref{tab:searchspace} specifies the complete tuning domain. The grids for
$\alpha$ and $\gamma$ are
$\{0.1+j(0.9/19):j=0,\ldots,19\}$; their retained values are displayed to six
decimal places in Table~\ref{tab:optuna_configurations}. The values
$T_F=0$ and $\kappa=0$ allow the full 2-opt and double-bridge stages,
respectively, to be disabled.

\begin{table}[H]
	\centering
	\small
	\setlength{\tabcolsep}{5pt}
	\caption{Optuna search space and study design.}
	\label{tab:searchspace}
	\begin{tabular}{ll}
		\toprule
		Parameter or control & Domain or fixed value \\
		\midrule
		$\alpha,\gamma$ &
		$\{0.1+j(0.9/19):j=0,\ldots,19\}$\\
		$P$ & $\{20,25,\ldots,60\}$\\
		$K$ & $\{10,15,\ldots,60\}$\\
		$s$ & $\{1,2,3,4\}$\\
		$\ell$ & $\{1,2,3,4,5\}$\\
		$T_2$ & $\{3,4,\ldots,12\}$\\
		$T_f$ & $\{20,50,100,200\}$\\
		$T_F$ & $\{0,100,500,1000\}$\\
		$\kappa$ & $\{0,5,10,15,20,25,30\}$\\
		$\mu$ & $\{500,1000,2000,3000,5000\}$\\
		\midrule
		$B$ (not optimized) & $100000$\\
		$\eta$ (not optimized) & $0.7$\\
		Tuning seeds (not optimized) & $\{42,123,999\}$\\
		Trials per instance & $50$\\
		Solver executions per instance & $150$\\
		\bottomrule
	\end{tabular}
\end{table}

Table~\ref{tab:optuna_configurations} reports the retained configurations and
the derived elite cardinality $E=\lceil\alpha P\rceil$. Including $E$ avoids
ambiguity between the real-valued elite fraction and the integer number of
elite particles used by the implementation. All retained configurations satisfy
$E<P$ and therefore initialize at least one particle randomly.

\begin{table}[H]
	\centering
	\caption{Instance-specific configurations retained by the Optuna studies.
		Values of $\alpha$ and $\gamma$ are displayed to six decimal places;
		$E=\lceil\alpha P\rceil$ is the effective elite cardinality computed from
		the full-precision value of $\alpha$.}
	\label{tab:optuna_configurations}
	\scalebox{0.70}{%
		\begin{tabular}{lrrrrrrrrrrrr}
			\toprule
			Instance & $\alpha$ & $E$ & $\gamma$ & $P$ & $K$ & $s$ & $\ell$ &
			$T_2$ & $T_f$ & $T_F$ & $\kappa$ & $\mu$\\
			\midrule
			d493
			& 0.905263 & 55 & 0.242105 & 60 & 55 & 2 & 1
			& 12 & 20 & 100 & 10 & 3000\\
			d657
			& 0.952632 & 24 & 0.336842 & 25 & 15 & 4 & 5
			& 8 & 20 & 100 & 25 & 3000\\
			rat783
			& 0.905263 & 55 & 0.289474 & 60 & 15 & 3 & 5
			& 5 & 20 & 0 & 5 & 2000\\
			pr1002
			& 0.952632 & 48 & 0.857895 & 50 & 30 & 2 & 3
			& 5 & 50 & 0 & 15 & 2000\\
			u1060
			& 0.715789 & 40 & 0.621053 & 55 & 40 & 1 & 2
			& 11 & 50 & 0 & 25 & 5000\\
			\bottomrule
		\end{tabular}%
	}
\end{table}

The variation in Table~\ref{tab:optuna_configurations} reflects the
instance-specific design of the studies. Because the parameters were optimized
jointly, their individual contributions are evaluated through the ablation
study in Section~\ref{sec:ablation}. The selection criterion prioritizes
solution quality, while runtime is measured independently in the final
experiments.

\subsection{Results and discussion}\label{sec:results}

Table~\ref{tab:performance_comparison} shows that ADPSO-ERLS attains the lowest
mean cost $\bar f$ and the lowest best cost $f_{\min}$ on all five instances.
Its relative error, measured as \eqref{eq:experimental-re}, ranges
from $3.04\%$ to $5.50\%$, while the best-solution Gap, measured as \eqref{eq:experimental-gap}, ranges from $1.93\%$ to $4.17\%$. The
closest competitor in best-solution quality is AACO-LST on d493, d657, rat783,
and pr1002, and DPSO6 on u1060; even so, the Gap margin separating
ADPSO-ERLS from this second-best method grows monotonically with instance size,
through the values $2.14$, $2.74$, $3.26$, $4.85$, and $5.62$ percentage points,
from d493 to u1060. The paired tests
reported below confirm that these cost differences are statistically significant
for every comparison.

For numerical consistency, every ADPSO-ERLS entry in the two main comparison
tables is computed from the same $50$ raw runs as the \textbf{full} row of
Table~\ref{tab:ablation-five-instances}. This variant activates all five
components, uses the instance-specific configuration in
Table~\ref{tab:optuna_configurations}, follows the common seed sequence, and
terminates at the strict global budget $B=100000$.

\begin{table}[H]
	\centering
	\caption{Performance comparison of ADPSO-ERLS with SA-DPSO, AACO-LST, DPSO6, DPSO-Q, and AAACO. The cost statistics $\bar f$, $f_{\min}$, and $f^\star$ follow the notation introduced with Eq.~\eqref{eq:sample-costs}--\eqref{eq:experimental-apd}. For each metric and instance, the best value is shown in bold and the second best is underlined.}
	\label{tab:performance_comparison}
	\scalebox{0.75}{
		\begin{tabular}{l|l|c|c|c|c|c|c|c}
			\toprule
			\textbf{Instance} & \textbf{Method} & $\boldsymbol{\bar f}$ & $\boldsymbol{f_{\min}}$ & $\boldsymbol{f^\star}$ & \textbf{Gap (\%)} & \textbf{RE (\%)} & \textbf{APD (\%)} & \textbf{Runtime (s)} \\
			\midrule
			\multirow{6}{*}{d493}
			& SA-DPSO &39707.10&38800.07&35002&10.85&13.44&\underline{2.34}&199.60\\
			& AACO-LST&\underline{37688.13}&\underline{36426.85}&35002&\underline{4.07}&\underline{7.67}&3.46&41.06\\
			& DPSO6&38514.00&37164.00&35002&6.18&10.03&3.63&\textbf{1.55}\\
			& DPSO-Q&40218.03&38922.24&35002&11.20&14.90&3.33&40.55\\
			& AAACO&43319.46&40915.96&35002&16.90&23.76&5.87&26.55\\
			& ADPSO-ERLS&\textbf{36067.72}&\textbf{35677.00}&35002&\textbf{1.93}&\textbf{3.04}&\textbf{1.10}&\underline{19.18}\\ \midrule
			\multirow{6}{*}{d657}
			& SA-DPSO&56788.29&55996.89&48912&14.48&16.10&\underline{1.41}&255.45\\
			& AACO-LST&\underline{53138.78}&\underline{51870.72}&48912&\underline{6.05}&\underline{8.64}&2.44&58.04\\
			& DPSO6&54674.98&53227.00&48912&8.82&11.78&2.72&\underline{2.15}\\
			& DPSO-Q&57385.05&56154.12&48912&14.81&17.32&2.19&98.32\\
			& AAACO&63001.63&60252.07&48912&23.18&28.81&4.56&46.31\\
			& ADPSO-ERLS&\textbf{51144.86}&\textbf{50531.00}&48912&\textbf{3.31}&\textbf{4.57}&\textbf{1.21}&\textbf{1.70}\\ \midrule
			\multirow{6}{*}{rat783}
			& SA-DPSO&9858.95&9770.56&8806&10.95&11.96&\textbf{0.90}&299.92\\
			& AACO-LST&\underline{9680.73}&\underline{9436.76}&8806&\underline{7.16}&\underline{9.93}&2.59&74.31\\
			& DPSO6&9944.48&9797.00&8806&11.25&12.93&1.51&\textbf{2.46}\\
			& DPSO-Q&10458.84&10238.00&8806&16.26&18.77&2.16&160.92\\
			& AAACO&11337.95&10764.78&8806&22.24&28.75&5.32&64.80\\
			& ADPSO-ERLS&\textbf{9240.22}&\textbf{9149.00}&8806&\textbf{3.90}&\textbf{4.93}&\underline{1.00}&\underline{2.51}\\ \midrule
			\multirow{6}{*}{pr1002}
			& SA-DPSO&316073.63&313634.11&259045&21.07&22.01&\textbf{0.78}&382.55\\
			& AACO-LST&292625.61&\underline{282283.07}&259045&\underline{8.97}&12.96&3.66&103.67\\
			& DPSO6&\underline{291319.14}&284086.00&259045&9.67&\underline{12.46}&2.55&\textbf{3.60}\\
			& DPSO-Q&305678.36&295875.12&259045&14.22&18.00&3.31&368.29\\
			& AAACO&346596.68&332098.32&259045&28.20&33.80&4.37&110.78\\
			& ADPSO-ERLS&\textbf{272269.58}&\textbf{269705.00}&259045&\textbf{4.12}&\textbf{5.11}&\underline{0.95}&\underline{12.51}\\ \midrule
			\multirow{6}{*}{u1060}
			& SA-DPSO&260655.21&258164.01&224094&15.20&16.32&\textbf{0.96}&428.15\\
			& AACO-LST&257559.85&250683.01&224094&11.87&14.93&2.74&114.54\\
			& DPSO6&\underline{252752.98}&\underline{246026.00}&224094&\underline{9.79}&\underline{12.79}&2.73&\textbf{3.27}\\
			& DPSO-Q&262901.90&256814.05&224094&14.60&17.32&2.37&413.70\\
			& AAACO&305607.95&291877.14&224094&30.25&36.37&4.70&116.07\\
			& ADPSO-ERLS&\textbf{236428.02}&\textbf{233446.00}&224094&\textbf{4.17}&\textbf{5.50}&\underline{1.28}&\underline{32.59}\\
			\bottomrule
	\end{tabular}}
\end{table}

The descriptive statistics in Table~\ref{tab:descriptive_statistics} reinforce
both the quality and the repeatability patterns. ADPSO-ERLS attains the smallest
standard deviation on d493 ($265.68$) and the second smallest on d657, rat783,
pr1002, and u1060, being outperformed in dispersion only by SA-DPSO on those
four instances. It is thus first or second in dispersion on every instance while
remaining the only method that simultaneously attains the smallest mean and best
cost throughout. Its APD lies between $0.95\%$ and $1.28\%$. SA-DPSO occasionally
records a smaller APD (e.g., rat783, pr1002, u1060), but around a materially
worse mean; low dispersion alone therefore does not imply proximity to the
optimum. Regarding runtime, DPSO6 is consistently the fastest method, whereas
ADPSO-ERLS trades additional computation for the lowest best and mean costs. On
d657 and rat783, ADPSO-ERLS requires only $1.70$ and $2.51$ seconds per run,
essentially matching DPSO6 ($2.15$ and $2.46$ seconds) while attaining
substantially lower costs. On d493, pr1002, and u1060, the more intensive tuned
refinement configurations raise the ADPSO-ERLS runtime to $19.18$, $12.51$, and
$32.59$ seconds, respectively; these values nevertheless remain below those of
SA-DPSO, DPSO-Q, and AAACO on the same instances, and also well below AACO-LST,
which requires $103.67$ and $114.54$ seconds on pr1002 and u1060 against
$12.51$ and $32.59$ seconds for ADPSO-ERLS, roughly $8.3$ and $3.5$ times faster.
The quality gain is thus not
obtained through uniformly higher runtime than every competitor: ADPSO-ERLS
occupies a favorable region of the quality--runtime trade-off, with DPSO6
defining the speed extreme. Since every method operates under $B=100000$ on
identical hardware, the runtime comparison is made under a common strict global
evaluation budget, and the remaining timing differences reflect the
computational cost of each method's operators and implementation.

\begin{table}[H]
	\centering
	\caption{Descriptive statistics obtained by ADPSO-ERLS and the comparison methods. The smallest reported standard deviation for each instance is shown in bold; the second smallest is underlined.}
	\label{tab:descriptive_statistics}
	\scalebox{0.75}{
		\begin{tabular}{l|l|c|c|c|c|c|c|c}
			\toprule
			\textbf{Instance} & \textbf{Method} & \textbf{Worst} & \textbf{Best} & \textbf{Mean} & \textbf{Median} & \textbf{Std} & \textbf{Lower95} & \textbf{Upper95} \\
			\midrule
			\multirow{6}{*}{d493}
			& SA-DPSO&40175.59&38800.07&39707.10&39818.37&\underline{351.45}&39609.68&39804.52\\
			& AACO-LST&38986.89&36426.85&37688.13&37670.33&576.88&37528.23&37848.03\\
			& DPSO6&40040.00&37164.00&38514.00&38543.00&633.24&38338.47&38689.53\\
			& DPSO-Q&41809.78&38922.24&40218.03&40196.95&640.53&40040.48&40395.58\\
			& AAACO&44915.31&40915.96&43319.46&43298.65&848.37&43084.30&43554.62\\
			& ADPSO-ERLS&36795.00&35677.00&36067.72&36028.50&\textbf{265.68}&35994.08&36141.36\\ \midrule
			\multirow{6}{*}{d657}
			& SA-DPSO&57360.93&55996.89&56788.29&56735.85&\textbf{296.63}&56706.07&56870.51\\
			& AACO-LST&54685.61&51870.72&53138.78&53149.70&703.93&52943.66&53333.90\\
			& DPSO6&55726.00&53227.00&54674.98&54664.00&543.42&54524.35&54825.61\\
			& DPSO-Q&59576.27&56154.12&57385.05&57330.76&731.27&57182.35&57587.75\\
			& AAACO&65990.84&60252.07&63001.63&63027.48&1381.94&62618.58&63384.68\\
			& ADPSO-ERLS&52548.00&50531.00&51144.86&51119.50&\underline{377.06}&51040.35&51249.37\\ \midrule
			\multirow{6}{*}{rat783}
			& SA-DPSO&9957.13&9770.56&9858.95&9847.70&\textbf{48.57}&9845.49&9872.41\\
			& AACO-LST&9973.51&9436.76&9680.73&9685.79&143.75&9640.88&9720.58\\
			& DPSO6&10164.00&9797.00&9944.48&9933.50&84.38&9921.09&9967.87\\
			& DPSO-Q&10650.95&10238.00&10458.84&10472.53&96.76&10432.02&10485.66\\
			& AAACO&11722.76&10764.78&11337.95&11367.66&214.16&11278.59&11397.31\\
			& ADPSO-ERLS&9358.00&9149.00&9240.22&9241.50&\underline{54.92}&9225.00&9255.44\\ \midrule
			\multirow{6}{*}{pr1002}
			& SA-DPSO&318358.12&313634.11&316073.63&315912.67&\textbf{1009.10}&315793.92&316353.34\\
			& AACO-LST&299384.07&282283.07&292625.61&292763.34&3849.49&291558.59&293692.63\\
			& DPSO6&296910.00&284086.00&291319.14&290867.50&3159.88&290443.27&292195.01\\
			& DPSO-Q&316033.62&295875.12&305678.36&305750.26&3647.07&304667.44&306689.28\\
			& AAACO&362611.29&332098.32&346596.68&346128.44&6270.30&344858.64&348334.72\\
			& ADPSO-ERLS&275220.00&269705.00&272269.58&272393.50&\underline{1628.81}&271818.11&272721.05\\ \midrule
			\multirow{6}{*}{u1060}
			& SA-DPSO&263656.05&258164.01&260655.21&260674.32&\textbf{1429.36}&260259.01&261051.41\\
			& AACO-LST&267280.32&250683.01&257559.85&257169.86&3507.90&256587.51&258532.19\\
			& DPSO6&258728.00&246026.00&252752.98&252819.00&2680.10&252010.09&253495.87\\
			& DPSO-Q&269281.99&256814.05&262901.90&262879.24&2725.92&262146.31&263657.49\\
			& AAACO&314250.46&291877.14&305607.95&306095.66&5013.51&304218.28&306997.62\\
			& ADPSO-ERLS&239825.00&233446.00&236428.02&236291.50&\underline{1447.13}&236026.90&236829.14\\
			\bottomrule
	\end{tabular}}
\end{table}

\subsection{Statistical analysis}\label{sec:statistical_analysis}

The analysis follows a global-test, post-hoc-comparison, and effect-size
sequence. The experimental unit for the omnibus comparison is a TSPLIB instance
($N=5$) and the treatment is an algorithm ($k=6$). The Friedman null hypothesis
states that all algorithms share the same rank distribution across instances.
Table~\ref{tab:friedman} reports $\chi_F^2(5)=23.057$ with
$p=3.292\times10^{-4}$, so the null hypothesis is rejected at $\alpha=0.05$.
Because the benchmark set is small, the less conservative Iman--Davenport
correction was also computed; it yields $F_F(5,20)=47.471$ with
$p=2.042\times10^{-10}$, confirming the global result. These omnibus tests
establish that at least one rank distribution differs, but neither identifies a
particular pair nor quantifies the magnitude of any cost difference.

\begin{table}[H]
	\centering
	\caption{Omnibus rank tests comparing the six optimization algorithms.}
	\label{tab:friedman}
	\renewcommand{\arraystretch}{1.2}
	\scalebox{0.75}{
	\begin{tabular}{lcccc}
		\toprule
		\textbf{Test} & \textbf{Statistic} & \textbf{Degrees of freedom} & \textbf{$p$-value} & \textbf{Decision} \\
		\midrule
		Friedman & $\chi_F^2=23.057$ & 5 & $3.292\times10^{-4}$ & Reject $H_0$ \\
		Iman--Davenport & $F_F=47.471$ & $(5,20)$ & $2.042\times10^{-10}$ & Reject $H_0$ \\
		\bottomrule
	\end{tabular}
}
\end{table}

The average ranks in Table~\ref{tab:meanrank} give the direction of the global
result. ADPSO-ERLS attains the lowest attainable mean rank ($1.0$), meaning it
ranks first on every benchmark entering the Friedman analysis. AACO-LST ($2.4$)
and DPSO6 ($2.8$) form the closest comparison group, while SA-DPSO ($4.0$) and
DPSO-Q ($4.8$) occupy intermediate positions. AAACO ($6.0$) obtains the least
favorable aggregate rank under these local implementations. Since a rank
encodes only relative ordering and not the numerical size of a cost difference,
it must be interpreted jointly with Gap,
RE, APD, and the descriptive statistics.

\begin{table}[H]
	\centering
	\caption{Average ranks obtained from the Friedman test (lower values indicate better performance).}
	\label{tab:meanrank}
	\renewcommand{\arraystretch}{1.2}
	\scalebox{0.75}{
	\begin{tabular}{lc}
		\hline
		\textbf{Algorithm} & \textbf{Mean Rank} \\
		\hline
		\textbf{ADPSO-ERLS} & \textbf{1.0} \\
		AACO-LST & 2.4 \\
		DPSO6 & 2.8 \\
		SA-DPSO & 4.0 \\
		DPSO-Q & 4.8 \\
		AAACO & 6.0 \\
		\hline
	\end{tabular}
}
\end{table}

The Friedman ranks identify ADPSO-ERLS as the first-ranked method on every
instance, but ranks do not measure practical separation. Accordingly,
Table~\ref{tab:wilcoxon_comparison} reports the $25$ paired Wilcoxon
signed-rank comparisons over the $50$ matched runs for each instance and
baseline. With paired differences oriented as baseline cost minus ADPSO-ERLS
cost, $24$ of the $25$ comparisons yield $W^-=0$, $W^+=1275$, and $Z=6.154$,
corresponding to $p=7.557\times10^{-10}$. The remaining comparison, ADPSO-ERLS
versus AACO-LST on d493, yields $W^-=1$, $W^+=1274$, $Z=6.144$, and
$p=8.031\times10^{-10}$. The largest raw $p$-value is therefore below
$8.1\times10^{-10}$, so even the Bonferroni bound for $25$ tests, at most
$2.008\times10^{-8}$, keeps every comparison significant under family-wise
error control at $\alpha=0.05$. ADPSO-ERLS achieves the lower cost in all $50$
matched runs in the $24$ tied comparisons, and in $49$ of the $50$ runs against
AACO-LST on d493.

To separate statistical significance from practical relevance, the pairwise
analysis also quantifies the magnitude of the observed differences. For each
comparison, the Hodges--Lehmann estimator gives the typical paired cost
difference, defined as ADPSO-ERLS minus the corresponding baseline, and its
$95\%$ confidence interval reflects the uncertainty of that estimate. All $25$
point estimates and all their confidence limits are negative, so every interval
excludes zero in the direction of lower cost for ADPSO-ERLS. Because TSPLIB
instances differ markedly in cost scale, these absolute reductions characterize
within-instance separation and are not compared directly across instances.

The distributional effect sizes support the same conclusion. Here the
Vargha--Delaney $A_{12}$ denotes the probability that an ADPSO-ERLS observation
exceeds a baseline observation, so lower values are favorable because tour cost
is minimized. In $23$ of the $25$ comparisons, $A_{12}=0$ and Cliff's
$\delta=-1$, indicating complete empirical separation in the observed samples.
The only nonzero values arise against AACO-LST, with
$(A_{12},\delta)=(0.0028,-0.9944)$ on d493 and $(0.0072,-0.9856)$ on d657; both
still represent near-complete separation and classify as large effects.
Consequently, the advantage of ADPSO-ERLS is supported not only by rejection of
the omnibus and pairwise null hypotheses, but also by uniformly negative
location shifts and large distributional effects across every evaluated
instance--algorithm pair.

\begin{table}[H]
	\centering
	\caption{Pairwise Wilcoxon signed-rank tests comparing ADPSO-ERLS with the competing algorithms over 50 independent runs.}
	\label{tab:wilcoxon_comparison}
	\renewcommand{\arraystretch}{1.15}
	\scalebox{0.75}{
		\begin{tabular}{lllrcccc}
			\toprule
			\textbf{Instance} &
			\textbf{Method A} &
			\textbf{Method B} &
			$\mathbf{W^-}$ &
			$\mathbf{W^+}$ &
			$\mathbf{Z}$ &
			$\mathbf{p}$-\textbf{value} &
			\textbf{Decision}\\
			\midrule
			
			\multirow{5}{*}{d493}
			& ADPSO-ERLS & SA-DPSO & 0 & 1275 & 6.154 & $7.557\times10^{-10}$ & Reject $H_0$\\
			& ADPSO-ERLS & AACO-LST & 1 & 1274 & 6.144 & $8.031\times10^{-10}$ & Reject $H_0$\\
			& ADPSO-ERLS & DPSO6 & 0 & 1275 & 6.154 & $7.555\times10^{-10}$ & Reject $H_0$\\
			& ADPSO-ERLS & DPSO-Q & 0 & 1275 & 6.154 & $7.557\times10^{-10}$ & Reject $H_0$\\
			& ADPSO-ERLS & AAACO & 0 & 1275 & 6.154 & $7.557\times10^{-10}$ & Reject $H_0$\\
			
			\midrule
			
			\multirow{5}{*}{d657}
			& ADPSO-ERLS & SA-DPSO & 0 & 1275 & 6.154 & $7.557\times10^{-10}$ & Reject $H_0$\\
			& ADPSO-ERLS & AACO-LST & 0 & 1275 & 6.154 & $7.557\times10^{-10}$ & Reject $H_0$\\
			& ADPSO-ERLS & DPSO6 & 0 & 1275 & 6.154 & $7.554\times10^{-10}$ & Reject $H_0$\\
			& ADPSO-ERLS & DPSO-Q & 0 & 1275 & 6.154 & $7.557\times10^{-10}$ & Reject $H_0$\\
			& ADPSO-ERLS & AAACO & 0 & 1275 & 6.154 & $7.557\times10^{-10}$ & Reject $H_0$\\
			
			\midrule
			
			\multirow{5}{*}{rat783}
			& ADPSO-ERLS & SA-DPSO & 0 & 1275 & 6.154 & $7.557\times10^{-10}$ & Reject $H_0$\\
			& ADPSO-ERLS & AACO-LST & 0 & 1275 & 6.154 & $7.557\times10^{-10}$ & Reject $H_0$\\
			& ADPSO-ERLS & DPSO6 & 0 & 1275 & 6.154 & $7.555\times10^{-10}$ & Reject $H_0$\\
			& ADPSO-ERLS & DPSO-Q & 0 & 1275 & 6.154 & $7.557\times10^{-10}$ & Reject $H_0$\\
			& ADPSO-ERLS & AAACO & 0 & 1275 & 6.154 & $7.557\times10^{-10}$ & Reject $H_0$\\
			
			\midrule
			
			\multirow{5}{*}{pr1002}
			& ADPSO-ERLS & SA-DPSO & 0 & 1275 & 6.154 & $7.557\times10^{-10}$ & Reject $H_0$\\
			& ADPSO-ERLS & AACO-LST & 0 & 1275 & 6.154 & $7.557\times10^{-10}$ & Reject $H_0$\\
			& ADPSO-ERLS & DPSO6 & 0 & 1275 & 6.154 & $7.557\times10^{-10}$ & Reject $H_0$\\
			& ADPSO-ERLS & DPSO-Q & 0 & 1275 & 6.154 & $7.557\times10^{-10}$ & Reject $H_0$\\
			& ADPSO-ERLS & AAACO & 0 & 1275 & 6.154 & $7.557\times10^{-10}$ & Reject $H_0$\\
			
			\midrule
			
			\multirow{5}{*}{u1060}
			& ADPSO-ERLS & SA-DPSO & 0 & 1275 & 6.154 & $7.557\times10^{-10}$ & Reject $H_0$\\
			& ADPSO-ERLS & AACO-LST & 0 & 1275 & 6.154 & $7.557\times10^{-10}$ & Reject $H_0$\\
			& ADPSO-ERLS & DPSO6 & 0 & 1275 & 6.154 & $7.557\times10^{-10}$ & Reject $H_0$\\
			& ADPSO-ERLS & DPSO-Q & 0 & 1275 & 6.154 & $7.557\times10^{-10}$ & Reject $H_0$\\
			& ADPSO-ERLS & AAACO & 0 & 1275 & 6.154 & $7.557\times10^{-10}$ & Reject $H_0$\\
			
			\bottomrule
	\end{tabular}}
\end{table}

\begin{table}[H]
	\centering
	\small
	\caption{Summary of the 25 pairwise effect-size analyses. Differences are oriented as ADPSO-ERLS cost minus baseline cost.}
	\label{tab:effect_size_summary}
	\scalebox{0.75}{
	\begin{tabular}{p{0.25\textwidth}p{0.65\textwidth}}
		\toprule
		\textbf{Measure} & \textbf{Result across the 25 comparisons} \\
		\midrule
		Hodges--Lehmann & All estimates and all 95\% confidence limits are negative. \\
		Vargha--Delaney $A_{12}$ & $0$ in 23 comparisons; $0.0028$ for d493 versus AACO-LST; $0.0072$ for d657 versus AACO-LST. \\
		Cliff's $\delta$ & $-1$ in 23 comparisons; $-0.9944$ and $-0.9856$ in the two remaining comparisons. \\
		Magnitude & Large in all 25 comparisons under the effect-size classification adopted in this study. \\
		Multiplicity check & All raw $p<8.1\times10^{-10}$; the largest Bonferroni-adjusted bound is below $2.008\times10^{-8}$. \\
		\bottomrule
	\end{tabular}
}
\end{table}

\subsection{Ablation study}\label{sec:ablation}

To determine whether the observed behavior depends on the coordinated
architecture rather than on a single operator, we conducted a
one-component-at-a-time ablation study on d493, d657, rat783, pr1002, and u1060.
The complete method was compared with five variants that respectively remove
heterogeneous initialization, evolutionary ERLS, final refinement,
double-bridge perturbation, or candidate-list restriction. Each variant
preserves the instance-specific configuration of
Table~\ref{tab:optuna_configurations}, uses the same 50 seeds, and terminates at
exactly $B=100000$ candidate-solution assessments. Distances follow the integer
TSPLIB \texttt{EUC\_2D} convention. The reported standard deviation uses the
sample estimator with denominator $n-1$, and runtime is the mean elapsed time
of the 50 individual runs rather than parallel wall time.

\begin{table}[H]
	\centering\small
	\caption{Descriptive ablation results on five TSPLIB instances (50 runs with common seeds, recorded limit $B=100000$). Lower values are better. $\Delta$Mean is relative to the complete method on the same instance. Runtime is mean elapsed time per run.}
	\label{tab:ablation-five-instances}
	\scalebox{0.65}{%
		\begin{tabular}{llrrrrrrrr}
			\toprule
			Instance & Variant & Best & Mean & Std & Gap (\%) & RE (\%) & APD (\%) & $\Delta$Mean (\%) & Runtime (s) \\
			\midrule
			d493 & \textbf{full} & 35677.00 & 36067.72 & 265.68 & 1.93 & 3.04 & 1.10 & --- & 19.18 \\
			& \texttt{no-heterogeneous-init} & 35941.00 & 36462.12 & 246.17 & 2.68 & 4.17 & 1.45 & +1.09 & 27.73 \\
			& \texttt{no-erls} & 36115.00 & 36644.66 & 284.42 & 3.18 & 4.69 & 1.47 & +1.60 & 0.27 \\
			& \texttt{no-final-refinement} & 35677.00 & 36444.08 & 418.77 & 1.93 & 4.12 & 2.15 & +1.04 & 21.39 \\
			& \texttt{no-double-bridge} & 35677.00 & 36256.02 & 422.31 & 1.93 & 3.58 & 1.62 & +0.52 & 21.28 \\
			& \texttt{no-candidate-lists} & 35579.00 & 36049.34 & 201.45 & 1.65 & 2.99 & 1.32 & -0.05 & 166.84 \\
			\midrule
			d657 & \textbf{full} & 50531.00 & 51144.86 & 377.06 & 3.31 & 4.57 & 1.21 & --- & 1.70 \\
			& \texttt{no-heterogeneous-init} & 51296.00 & 52143.82 & 499.64 & 4.87 & 6.61 & 1.65 & +1.95 & 0.82 \\
			& \texttt{no-erls} & 50455.00 & 51412.54 & 359.82 & 3.15 & 5.11 & 1.90 & +0.52 & 0.35 \\
			& \texttt{no-final-refinement} & 50913.00 & 52309.24 & 1031.40 & 4.09 & 6.95 & 2.74 & +2.28 & 1.60 \\
			& \texttt{no-double-bridge} & 50772.00 & 51414.46 & 447.23 & 3.80 & 5.12 & 1.27 & +0.53 & 1.61 \\
			& \texttt{no-candidate-lists} & 50497.00 & 51331.16 & 302.62 & 3.24 & 4.95 & 1.65 & +0.36 & 54.47 \\
			\midrule
			rat783 & \textbf{full} & 9149.00 & 9240.22 & 54.92 & 3.90 & 4.93 & 1.00 & --- & 2.51 \\
			& \texttt{no-heterogeneous-init} & 9273.00 & 9483.84 & 74.97 & 5.30 & 7.70 & 2.27 & +2.64 & 0.52 \\
			& \texttt{no-erls} & 9217.00 & 9311.94 & 49.69 & 4.67 & 5.75 & 1.03 & +0.78 & 0.33 \\
			& \texttt{no-final-refinement} & 9159.00 & 9327.88 & 81.03 & 4.01 & 5.93 & 1.84 & +0.95 & 2.46 \\
			& \texttt{no-double-bridge} & 9149.00 & 9249.62 & 53.98 & 3.90 & 5.04 & 1.10 & +0.10 & 2.38 \\
			& \texttt{no-candidate-lists} & 9173.00 & 9301.06 & 51.95 & 4.17 & 5.62 & 1.40 & +0.66 & 94.98 \\
			\midrule
			pr1002 & \textbf{full} & 269705.00 & 272269.58 & 1628.81 & 4.12 & 5.11 & 0.95 & --- & 12.51 \\
			& \texttt{no-heterogeneous-init} & 274978.00 & 279939.54 & 2601.34 & 6.15 & 8.07 & 1.80 & +2.82 & 1.95 \\
			& \texttt{no-erls} & 270795.00 & 274807.84 & 2068.68 & 4.54 & 6.08 & 1.48 & +0.93 & 0.74 \\
			& \texttt{no-final-refinement} & 270969.00 & 274519.90 & 1856.13 & 4.60 & 5.97 & 1.31 & +0.83 & 12.30 \\
			& \texttt{no-double-bridge} & 269712.00 & 272786.56 & 1800.19 & 4.12 & 5.30 & 1.14 & +0.19 & 12.31 \\
			& \texttt{no-candidate-lists} & 267934.00 & 271510.90 & 1746.59 & 3.43 & 4.81 & 1.33 & -0.28 & 342.88 \\
			\midrule
			u1060 & \textbf{full} & 233446.00 & 236428.02 & 1447.13 & 4.17 & 5.50 & 1.28 & --- & 32.59 \\
			& \texttt{no-heterogeneous-init} & 235356.00 & 239368.60 & 1754.95 & 5.03 & 6.82 & 1.70 & +1.24 & 43.15 \\
			& \texttt{no-erls} & 235880.00 & 238527.76 & 1429.86 & 5.26 & 6.44 & 1.12 & +0.89 & 1.36 \\
			& \texttt{no-final-refinement} & 237069.00 & 242286.92 & 3125.48 & 5.79 & 8.12 & 2.20 & +2.48 & 31.88 \\
			& \texttt{no-double-bridge} & 233880.00 & 238345.66 & 2118.61 & 4.37 & 6.36 & 1.91 & +0.81 & 31.88 \\
			& \texttt{no-candidate-lists} & 232731.00 & 235325.02 & 1515.88 & 3.85 & 5.01 & 1.11 & -0.47 & 740.24 \\
			\bottomrule
		\end{tabular}%
	}
\end{table}

Removing heterogeneous initialization increases mean cost on every instance,
from 1.09\% on d493 to 2.82\% on pr1002. This consistent association supports,
but does not by itself prove, a beneficial role for combining constructive and
random initial tours under the tested configurations. Disabling the
evolutionary ERLS stage also worsens every mean, by 0.52--1.60\%, while removing
the final refinement produces degradations of 0.83--2.48\% and the largest
dispersion on d657 and u1060. This pattern is consistent with a contribution
of incumbent-focused postprocessing to expected quality and repeatability, but
component interactions cannot be separated by a one-at-a-time design.
Suppressing double-bridge perturbation has a smaller but consistent effect of
0.10--0.81\%, suggesting a complementary rather than dominant association with
the attained quality.

Candidate lists exhibit a different trade-off. Replacing them with the complete
neighborhood changes mean cost by only $-0.47$\% to $+0.66$\%, and slightly
improves the mean on d493, pr1002, and u1060. However, average runtime rises from
19.18 to 166.84 seconds on d493, from 1.70 to 54.47 seconds on d657, from 2.51
to 94.98 seconds on rat783, from 12.51 to 342.88 seconds on pr1002, and from
32.59 to 740.24 seconds on u1060. Candidate restriction should therefore be
understood primarily as a computational-efficiency mechanism: it preserves
competitive solution quality while avoiding the rapidly increasing work of a
complete 2-opt neighborhood. This contrast also shows why equal FES alone does
not equalize all internal search work and why runtime must accompany the
quality metrics.

Overall, no single ablated component explains the complete pattern. The three
quality-oriented mechanisms have consistent but different magnitudes, whereas
candidate lists control the cost of obtaining those improvements. The evidence
therefore supports the architectural interpretation of ADPSO-ERLS as a
coordination of complementary exploration, intensification, refinement, and
neighborhood-restriction mechanisms.

\subsection{Scalability}
The scalability results cover four TSPLIB instances, namely
fl1400, u1432, u2319, and usa13509
\cite{reinelt1991}, together with three instances from the Waterloo National
TSP collection, namely Canada4663, Japan9847, and
Italy16862 \cite{waterlooNationalTSP}. Each row in
Table~\ref{tab:scalability_results} summarizes 50 independent runs. For the
TSPLIB instances, $f^\star$ denotes the reference optimum reported by TSPLIB,
whereas for the Waterloo instances it denotes the optimal tour value published
for the corresponding national instance. All experiments used the same
population-stage stopping threshold, $B=100000$, as in the main comparison.
This common threshold is reported explicitly for every instance in
Table~\ref{tab:scalability_results}.

\begin{table}[H]
	\centering
	\caption{ADPSO-ERLS scalability results from 50 independent runs per instance under the fixed limit $B=100000$.}
	\label{tab:scalability_results}
	\scalebox{0.75}{\begin{tabular}{lrrrrrrr}
			\toprule
			Instance&Cities&$B$&$f_{\min}$&$\bar f$&Gap (\%)&RE (\%)&Runtime/run (s)\\\midrule
			fl1400&1400&100000&20555&20677.50&2.13&2.74&31.64\\
			u1432&1432&100000&160478&161480.40&4.91&5.56&31.08\\
			u2319&2319&100000&242043&243081.46&3.32&3.77&6.52\\
			Canada4663&4663&100000&1356654&1378488.88&5.14&6.83&242.70\\
			Japan9847&9847&100000&523473&527296.88&6.41&7.19&541.44\\
			usa13509&13509&100000&21191580&21310314.20&6.05&6.64&600.59\\
			Italy16862&16862&100000&594163&597137.36&6.61&7.15&659.33\\
			\bottomrule\end{tabular}}
\end{table}

For the three largest instances, ADPSO-ERLS attains Gap values between 6.05\%
and 6.61\%, with mean runtimes between 541 and 659 s per run. In particular,
the method processes the 16862-city instance Italy16862 in less than
11 minutes on the reported computational platform while maintaining a Gap
below 6.7\%. These results extend the empirical evaluation beyond the
medium-sized TSPLIB instances used in the comparative study and demonstrate
that the proposed architecture remains computationally viable on substantially
larger search spaces.

The 6.52 s recorded for u2319 also indicates that the observed runtime
is determined by the interaction among instance structure, population size,
candidate-list size, and local-refinement settings, rather than by the number
of cities alone. All runs use the common strict global budget
$B=100000$, while the remaining controls are instance-specific. The runtime
results are therefore interpreted as direct measurements of practical
performance under the reported configurations, rather than as samples from a
single controlled scaling experiment.

This empirical behavior is consistent with the algorithmic design described in
Section~\ref{sec:algorithm}. Candidate lists restrict the neighborhoods explored
by the principal local-search stages, while explicit pass and improvement
limits bound the amount of refinement applied during a run. The FES accounting
of Section~\ref{sec:fes-accounting} continues to apply without modification on
these instances: a candidate-restricted 2-opt pass inspects at most $O(mK)$
moves and charges one FES per inspected move, so for the largest instances a
single full restricted pass would already exceed the global budget. Consequently,
the scalability runs use small instance-specific candidate-list sizes and rely on
the interruption rule of Eq.~\eqref{eq:fes_invariant}: the budgeted local-search
stages perform only partial neighborhood scans and stop cleanly at their active
deadline, committing every improvement accepted before interruption. The reported
Gaps are therefore produced by heterogeneous initialization followed by partial
budgeted refinement under the same strict limit $B=100000$, not by exhaustive
neighborhood exploration. The main structural
costs remain the $O(m^2)$ memory required by the stored distance matrix and the
$O(m^2)$ worst-case work of an unrestricted 2-opt pass. Within these
constraints, the experiments provide evidence that ADPSO-ERLS can produce
competitive-quality tours for instances ranging from 1400 to 16862 cities on
the stated hardware. This conclusion concerns practical feasibility over the
tested range and does not require extrapolation through an empirical
asymptotic model.

\begin{figure}[H]
\centering
\begin{subfigure}{0.48\textwidth}\includegraphics[width=\linewidth]{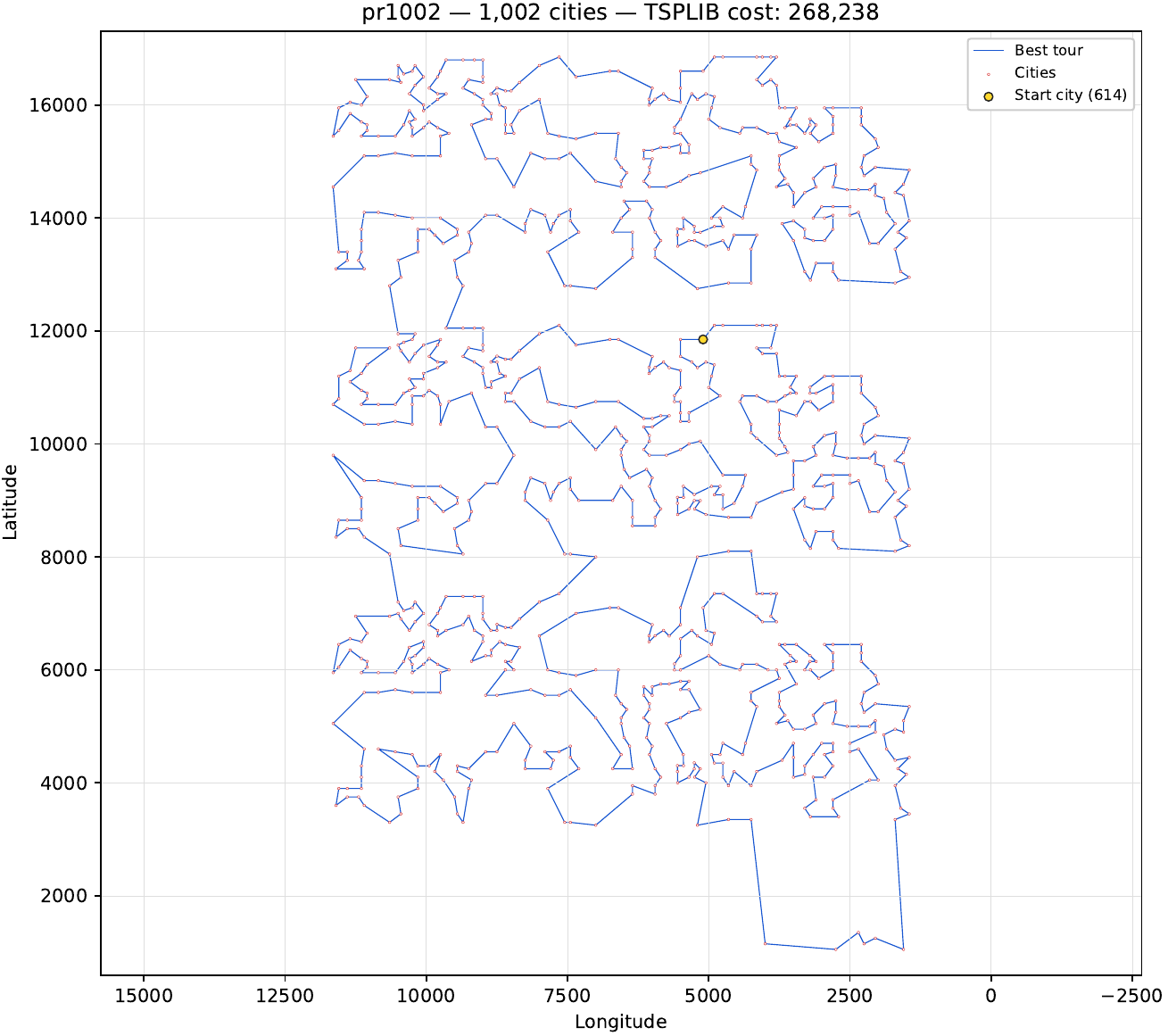}\caption{pr1002}\end{subfigure}\hfill
\begin{subfigure}{0.48\textwidth}\includegraphics[width=\linewidth]{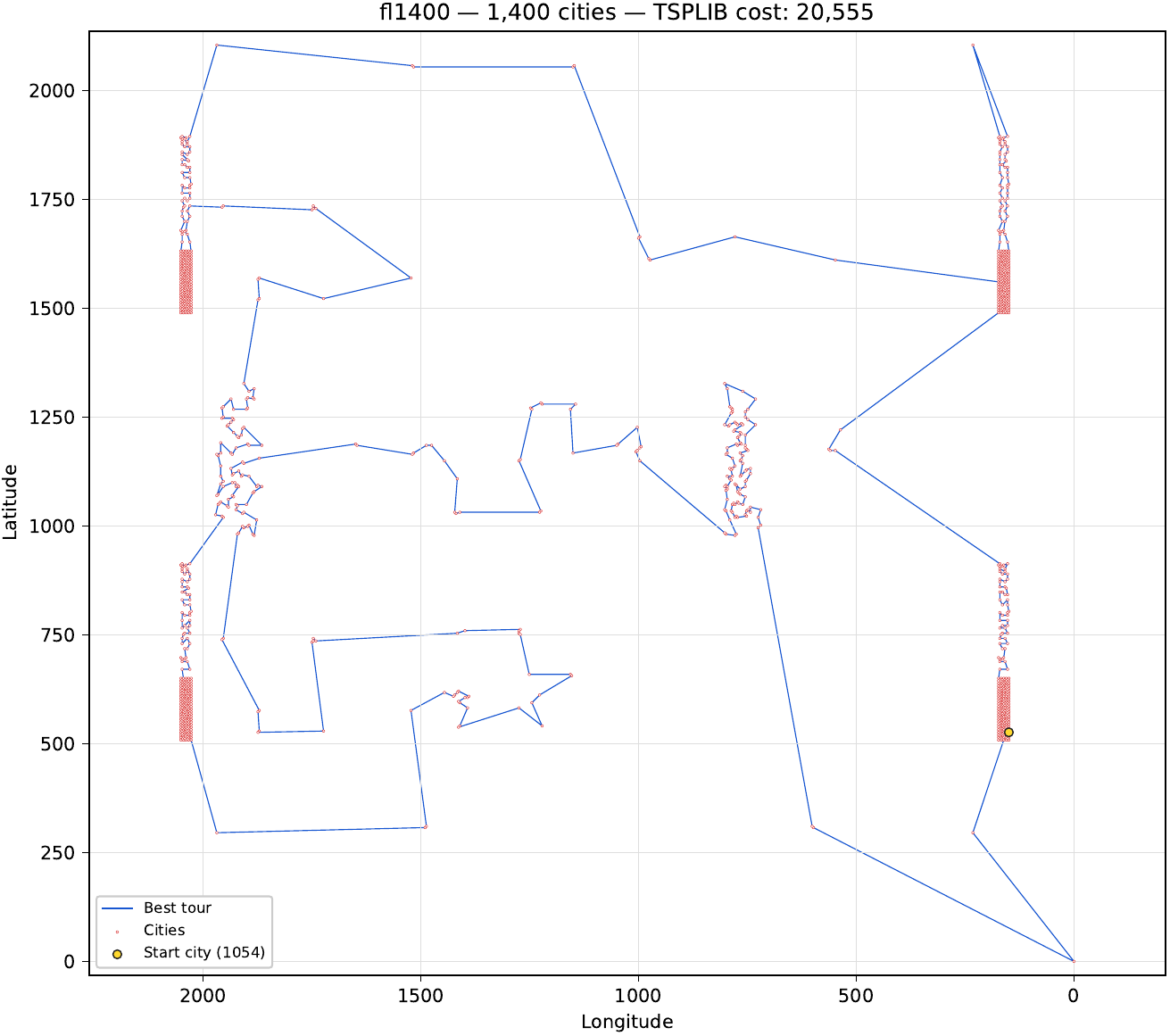}\caption{fl1400}\end{subfigure}
\begin{subfigure}{0.48\textwidth}\includegraphics[width=\linewidth]{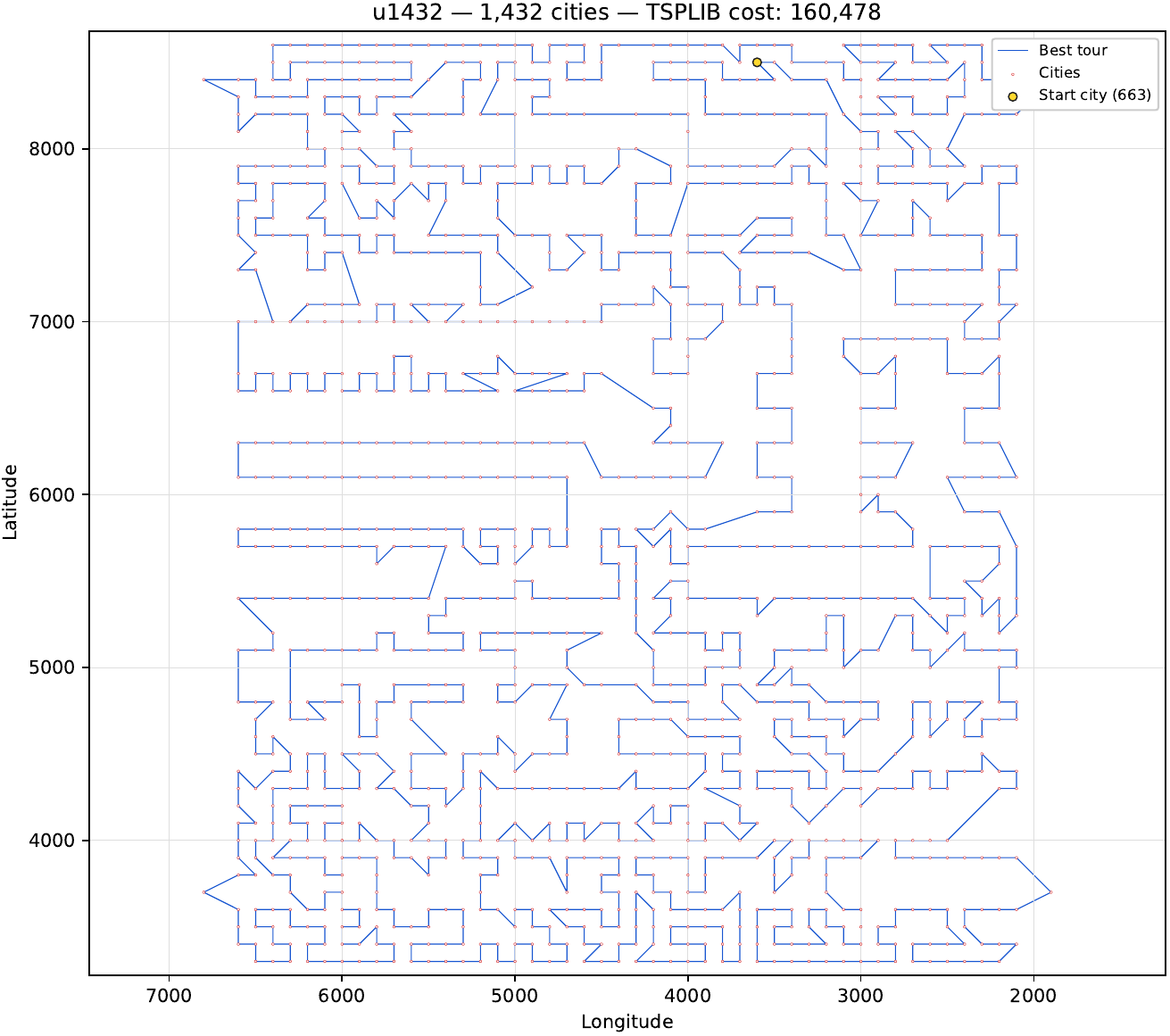}\caption{u1432}\end{subfigure}\hfill
\begin{subfigure}{0.48\textwidth}\includegraphics[width=\linewidth]{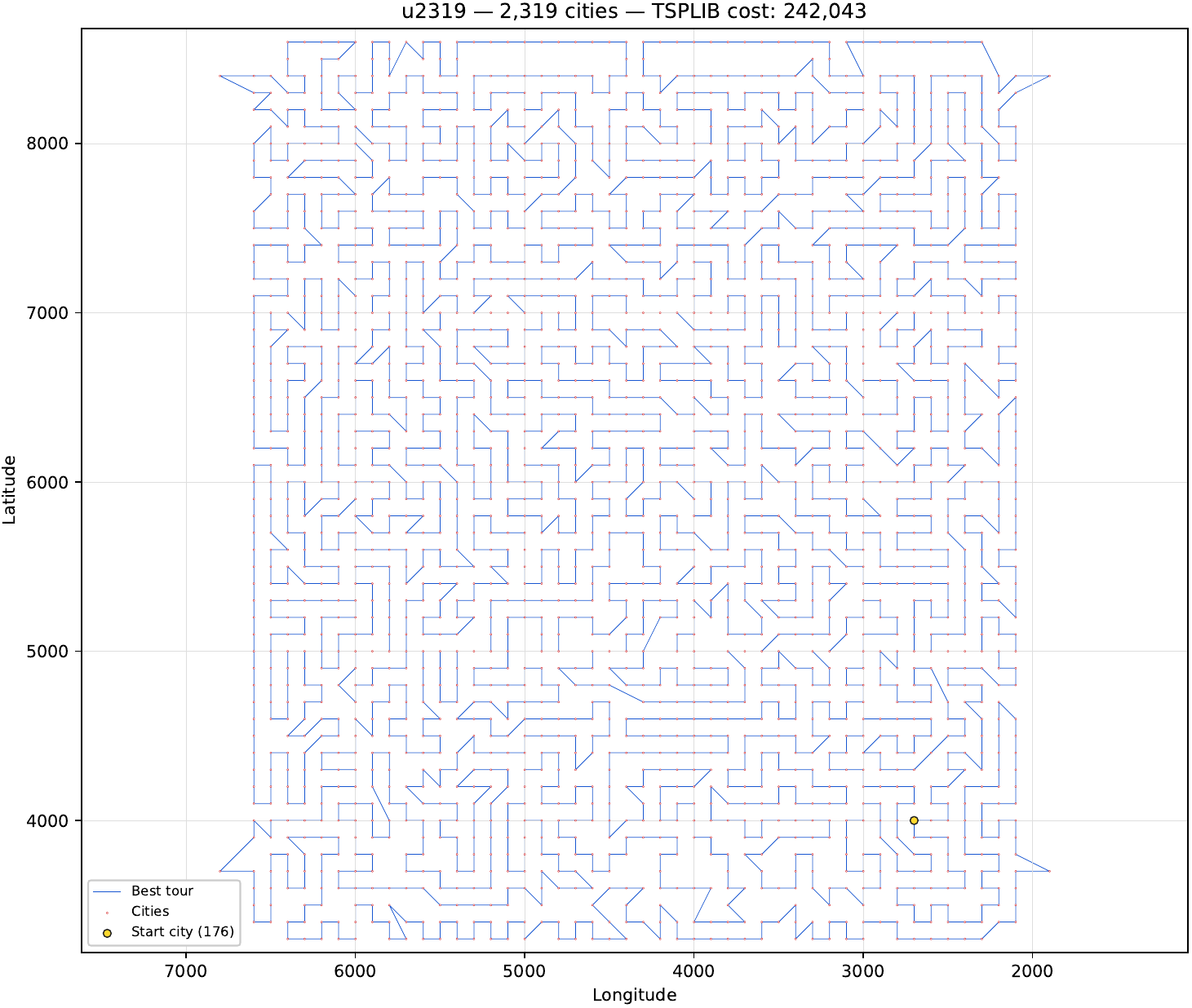}\caption{u2319}\end{subfigure}
\caption{ADPSO-ERLS tours for planar instances with 1,002--2,319
cities. Axes denote the source coordinates $x$ and $y$.}\label{fig:scalability_tours_a}
\end{figure}

\begin{figure}[H]
\centering
\begin{subfigure}{0.48\textwidth}\includegraphics[width=\linewidth]{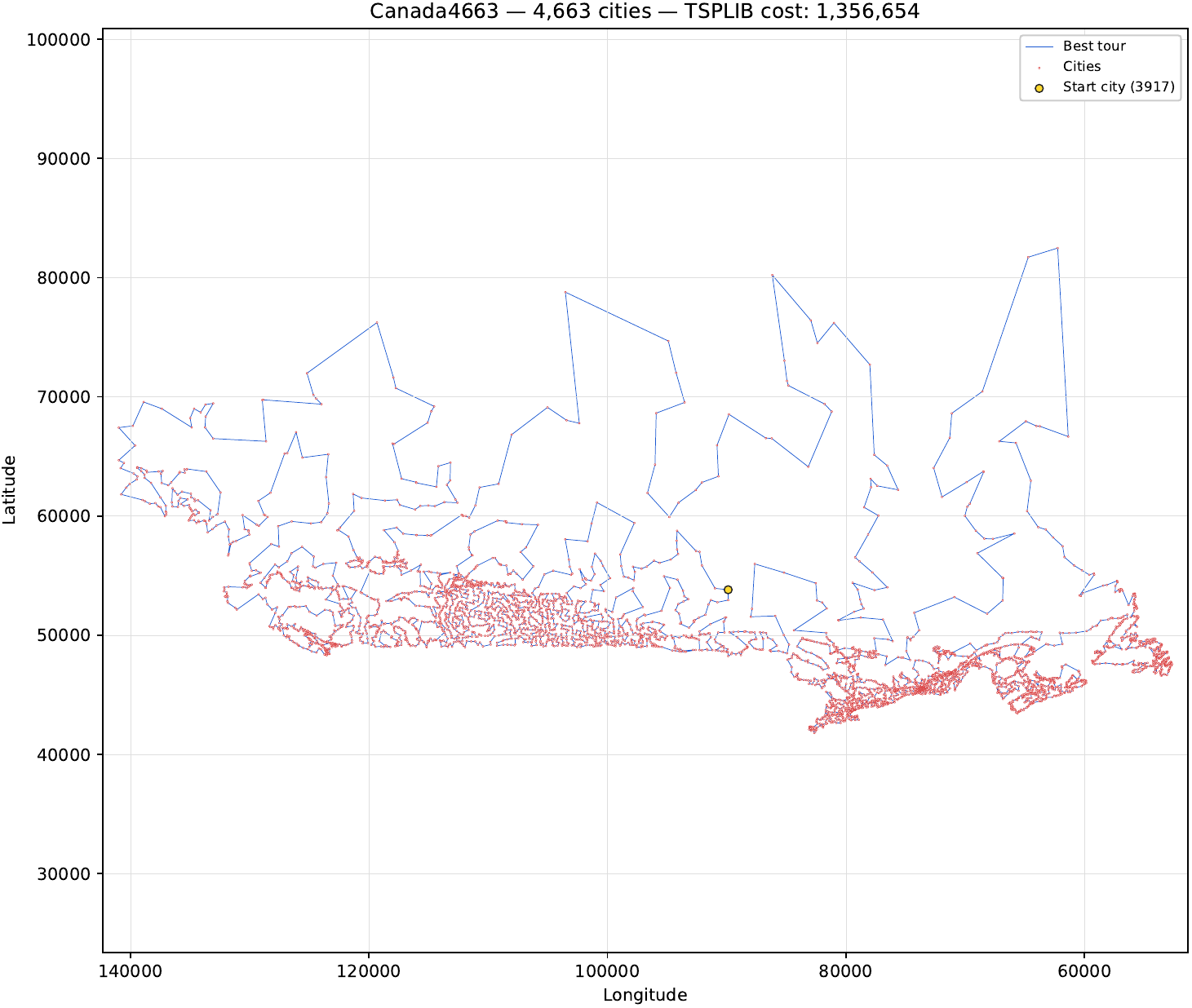}\caption{Canada4663}\end{subfigure}\hfill
\begin{subfigure}{0.48\textwidth}\includegraphics[width=\linewidth]{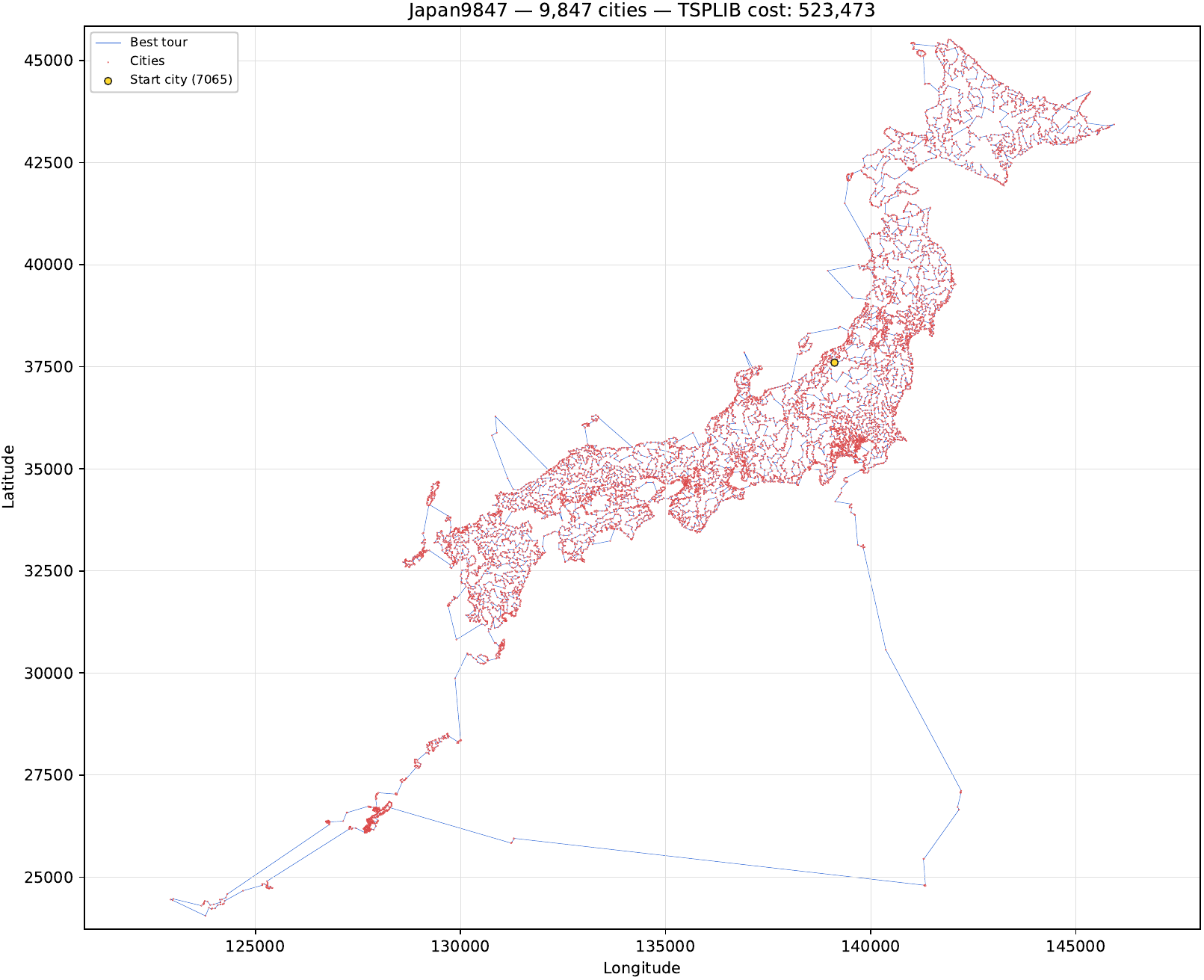}\caption{Japan9847}\end{subfigure}
\begin{subfigure}{0.48\textwidth}\includegraphics[width=\linewidth]{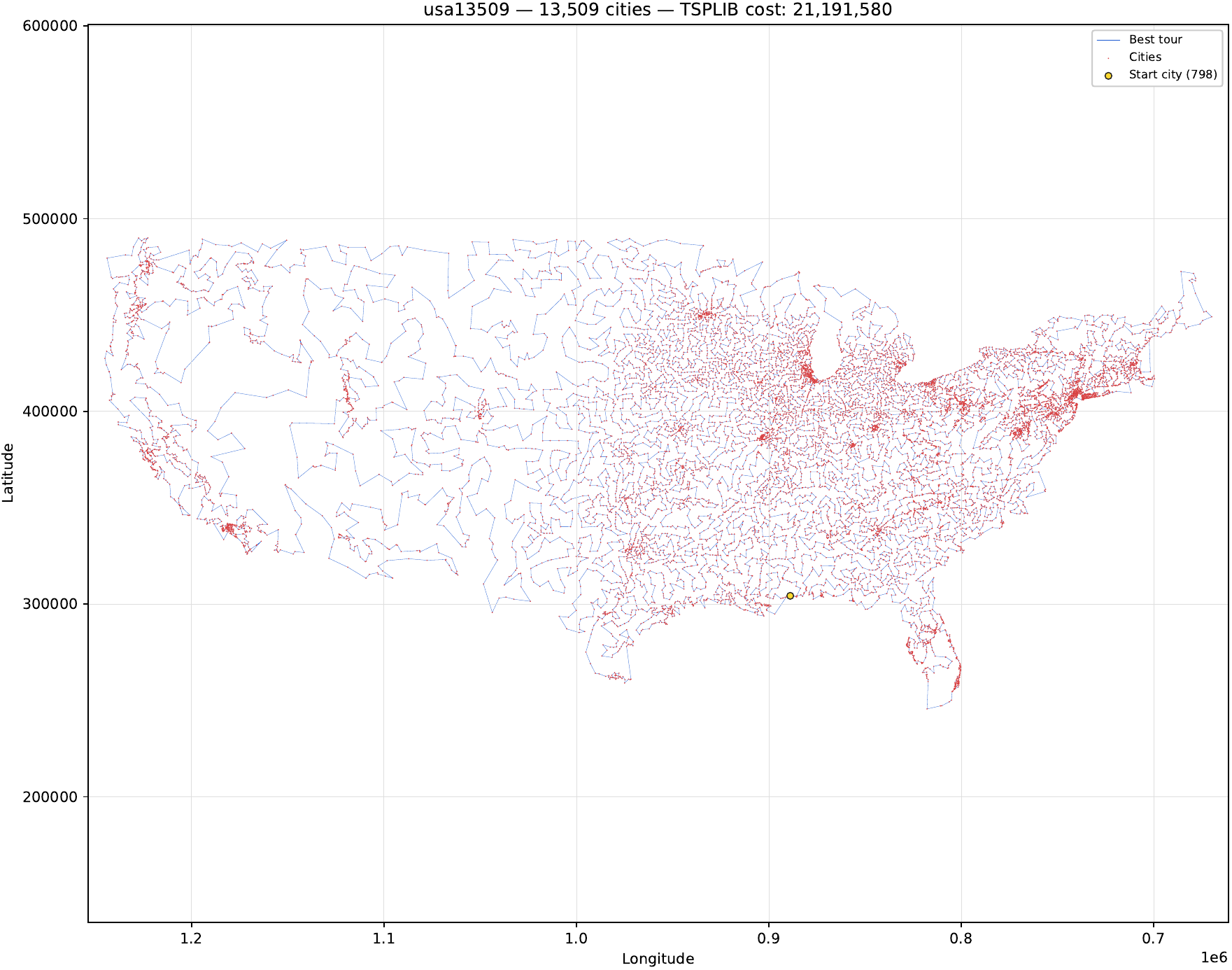}\caption{usa13509}\end{subfigure}\hfill
\begin{subfigure}{0.48\textwidth}\includegraphics[width=\linewidth]{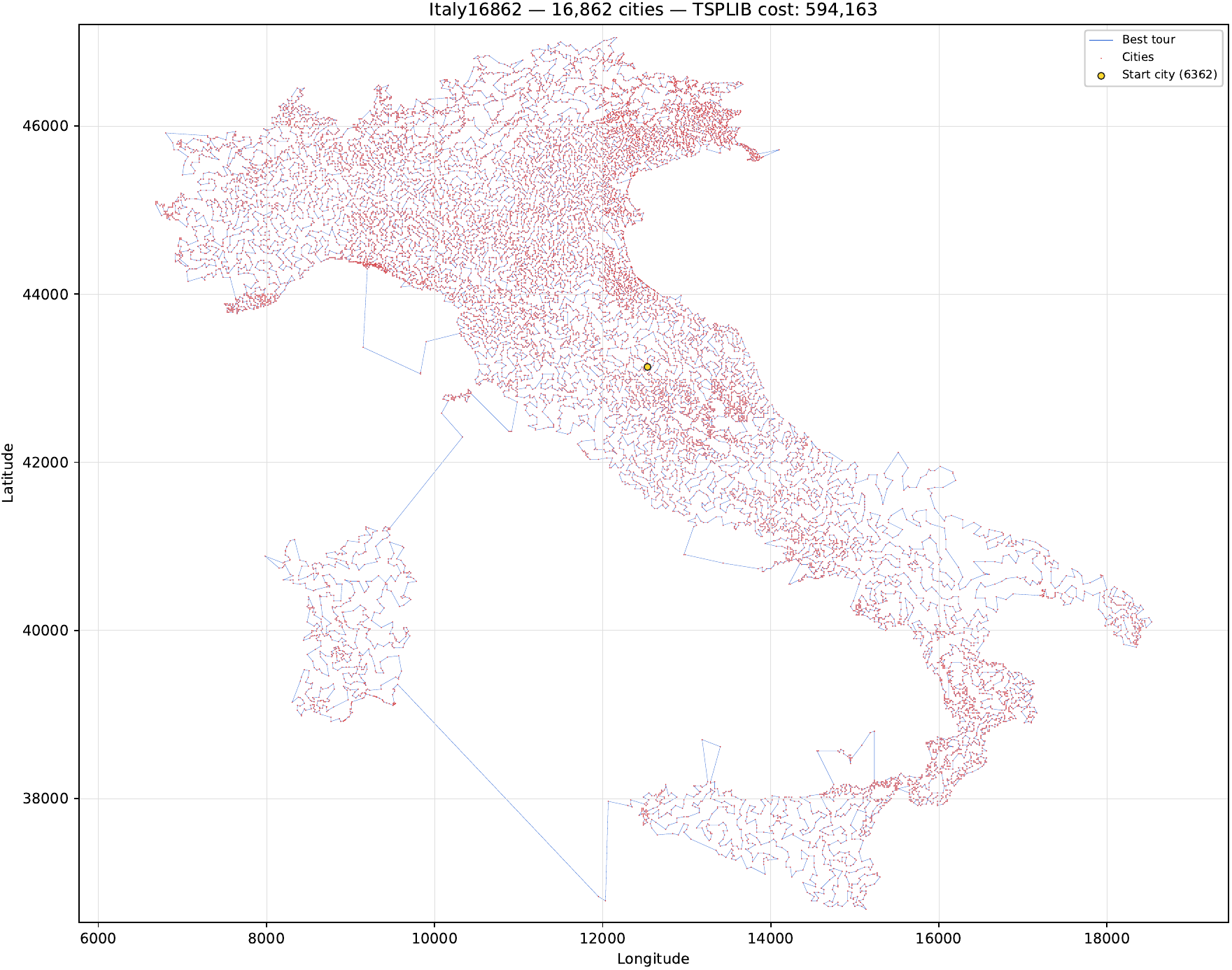}\caption{Italy16862}\end{subfigure}
\caption{ADPSO-ERLS tours for larger instances. Canada4663,
Japan9847, and Italy16862 use geographic axes; usa13509 retains the planar
TSPLIB labels $x$ and $y$.}\label{fig:scalability_tours_b}
\end{figure}

Figures~\ref{fig:scalability_tours_a} and
\ref{fig:scalability_tours_b} complement the numerical results by displaying
eight tours, including pr1002 as a planar reference. TSPLIB
instances are represented using their planar coordinates $(x,y)$, whereas the
Waterloo national instances use the geographic coordinates supplied with the
source data. The orientation of each axis follows the corresponding source
convention and does not affect the objective value.


\section{Conclusions}\label{sec:conclusions}

This study addressed a design question that is often left implicit in hybrid
discrete swarm optimization for the Traveling Salesman Problem: under a fixed
computational budget, how should search effort be distributed between
population-level exploration and progressively stronger local intensification?
Rather than adding another movement rule or neighborhood operator, ADPSO-ERLS
answers this question architecturally, exposing the allocation of effort as a set
of independently configurable parameters over memory-guided permutation search,
candidate-restricted elitist local improvement, and an incumbent-only final
stage combining candidate-restricted 2-opt, optional full 2-opt, and
double-bridge perturbations.

The experimental evidence consistently supports this design. Under the common
strict limit of $B=100000$ candidate-solution assessments, with all six
algorithms implemented in Rust and executed on identical hardware, ADPSO-ERLS
attained the lowest best and lowest mean cost on every TSPLIB instance
considered. Best-tour Gap values of $1.93$--$4.17\%$ and relative-error values of
$3.04$--$5.50\%$ show that the advantage holds for both the best tour and the
expected performance over repeated runs, rather than depending on isolated
favorable executions. The method also displayed the smallest standard deviation
on d493 and the second smallest on the remaining four instances, so improved
quality was accompanied by competitive run-to-run stability, and its APD stayed
between $0.95\%$ and $1.28\%$. The Friedman test assigned ADPSO-ERLS the best
attainable mean rank of $1.0$, the Iman--Davenport correction confirmed the
omnibus rejection, and all twenty-five paired Wilcoxon comparisons favored the
proposed method after family-wise error control. Every Hodges--Lehmann interval
was negative, and the Vargha--Delaney and Cliff effect sizes indicated large,
near-complete empirical separation from the baselines. Because the comparison was
conducted under a shared language, hardware, and stopping rule, these gains
reflect the architecture rather than a larger evaluation budget: ADPSO-ERLS was
faster than four of the five competitors on every instance and fastest of all on
d657, with DPSO6 defining the speed extreme.

The ablation study reinforces this interpretation. Removing heterogeneous
initialization raised mean cost by $1.09$--$2.82\%$, disabling evolutionary ERLS
by $0.52$--$1.60\%$, removing final refinement by $0.83$--$2.48\%$, and
suppressing double-bridge perturbation by $0.10$--$0.81\%$. Candidate lists acted
mainly as an efficiency mechanism: complete neighborhoods changed mean quality by
less than $0.7\%$ in either direction while inflating average runtime by factors
of roughly $8.7$ on d493 up to $37.8$ on rat783. No single ablated component
explains the full pattern, which supports the reading of ADPSO-ERLS as a
coordination of complementary mechanisms rather than a dependence on one operator.

The scalability experiments extend these observations beyond the primary
benchmark set. Using the same strict budget, ADPSO-ERLS was applied to four
TSPLIB and three Waterloo National TSP instances of $1{,}400$ to $16{,}862$
cities while keeping best-tour Gaps below $6.7\%$, processing the
$16{,}862$-city Italy16862 instance in under eleven minutes on the reported
platform. These experiments demonstrate practical viability across substantially
larger dimensions; they are not intended to establish asymptotic scaling laws.

Beyond the numerical improvements, the principal contribution is methodological.
The individual ingredients, namely discrete swarm learning, candidate-list
restriction, local search, perturbation, and hybridization, are established, and
no novelty is claimed for them in isolation. Instead, treating the allocation of
search effort itself as an explicit, auditable, and optimizable design variable
yields a transparent and reproducible framework, and suggests that future
advances may depend less on new operators than on principled mechanisms that
decide where, when, how often, and how strongly existing operators consume the
budget.

\section*{Acknowledgements}

The authors would like to thank the Department of Quantitative Methods at Universidad Loyola Andalusia for providing an excellent research environment and institutional support during the development of this work.

\section*{Declarations}

\subsection*{Ethical Approval}
Not applicable.

\subsection*{Consent to Participate}
Not applicable.

\subsection*{Consent to Publish}
Not applicable.

\subsection*{Data Availability Statement}
The experimental archive supporting the findings of this study, comprising the
per-run logs of the 50 independent executions per instance, the recorded budget
traces, the paired Wilcoxon comparisons, the archived best tours, and the tuning
driver, is available from the corresponding author upon reasonable request. All
benchmark instances are taken from the public TSPLIB and Waterloo National TSP
libraries cited in the text.

\subsection*{Author Contributions}
Conceptualization, A.S.-L. and J.-M.; methodology, A.S.-L., J.-M. and L.A. -B., software, A.S.-L., J.-M. and JM. -V.; validation, A.S.-L., J.-M. and JM. -V.; formal analysis, A.S.-L. and J.-M.; investigation, A.S.-L., J.-M. and JM. -V.; resources, A.S.-L.; data curation, A.S.-L., J.-M., and JM. -V. and L.A. -B.; writing--original draft preparation, A.S.-L. and JM. -V.; writing--review and editing, J.-M.; visualization, A.S.-L. and J.M; supervision, A.S.-L., J. M. and JM. -V.; project administration, A.S.-L. All authors have read and agreed to the published version of the manuscript.

\subsection*{Funding}
This research received no external funding.

\subsection*{Competing Interests}
The authors declare that they have no competing interests.


\end{document}